\documentclass[11pt,a4paper]{article}

\usepackage{graphicx}  
\usepackage{amsmath,amsfonts,amssymb}  
\usepackage{color}
\usepackage{url}
\usepackage{mathrsfs}
\usepackage{algorithmic}
\usepackage[ruled,vlined]{algorithm2e}

\newtheorem{definition}{Definition}
\newtheorem{conjecture}{Conjecture}

\newcommand{\bbN}{\mathbb{N}}
\newcommand{\bbZ}{\mathbb{Z}}
\newcommand{\bbR}{\mathbb{R}}
\newcommand{\bbC}{\mathbb{C}}
\newcommand{\bbS}{\mathbb{S}}

\newcommand{\n}{\mathbf{n}}

\begin{document}

\begin{center}
{\Huge Computing planar shape and critical point evolution under curvature-driven flows}

\vspace*{5mm}

{\Large Eszter Feh{\'e}r$^{1,2}$, Gábor Domokos$\,^{1,2}$, and Bernd Krauskopf$^3$}
\end{center}

\vspace*{1mm}

\noindent
{$^1$ MTA-BME Morphodynamics  Research Group \\[1mm] 
$^2$ Department of Mechanics, Materials and Structures, Budapest University of Technology and Economics, M{\H{u}}egyetem rakpart 1-3. K.II.61., 1111 Budapest, Hungary \\[1mm]
$^3$ Department of Mathematics, University of
  Auckland, Private Bag 92019, Auckland 1142, New Zealand}

\vspace*{2mm}

\begin{center}
{\large October 2020}
\end{center}

\vspace*{3mm}

\begin{abstract}
We are concerned with the evolution of planar, star-like curves and associated shapes under a broad class of curvature-driven geometric flows, which we refer to as the Andrews-Bloore flows. This family of flows has two parameters that control one constant and one curvature-dependent component for the velocity in the direction of the normal to the curve. The Andrews-Bloore flow includes as special cases the well known Eikonal, curve-shortening and affine shortening flows, and for positive parameter values its evolution shrinks the area enclosed by the curve to zero in finite time. A question of key interest has been how various shape descriptors of the evolving shape behave as this limit is approached. Star-like curves (which include convex curves) can be represented by a periodic scalar polar distance function $r(\varphi)$ measured from a reference point, which may or may not be fixed. An important question is how the numbers and the trajectories of critical points of the distance function $r(\varphi)$ and of the curvature $\kappa(\varphi)$ (characterized by $dr/d\varphi=0$ and $d\kappa /d\varphi=0$, respectively) evolve under the Andrews-Bloore flows for different choices of the parameters. \\ 
We present a numerical method that is specifically designed to meet the challenge of computing accurate trajectories of the critical points of an evolving curve up to the vicinity of a limiting shape. Each curve is represented by a piecewise polynomial periodic distance function, as determined by a chosen mesh; different types of meshes and mesh adaptation can be chosen to ensure a good balance between accuracy and computational cost. As we demonstrate with benchmark tests and two longer case studies, our method allows one to perform numerical investigations into subtle questions of planar curve evolution. More specifically --- in the spirit of experimental mathematics --- we provide illustrations of some known results, numerical evidence for two stated conjectures, as well as new insights and observations regarding the limits of shapes and their critical points.
\end{abstract}

\section{Introduction}

The curve-shortening flow describes the propagation of an embedded, closed planar curve $\Gamma$ such that individual points move in the direction of the inward normal with speed $v$ proportional to the curvature $\kappa$:
\begin{equation}
\label{eq:curveshortening}
v=C\kappa.
\end{equation}
The curvature-shortening flow \ref{eq:curveshortening} belongs to a broad class of nonlinear partial differential equations (PDEs) called \emph{curvature-driven flows} where the speed of evolution in the normal direction is given as some function of the curvature (in two dimensions) or curvatures (in higher dimensions). While locally defined, curvature-driven flows have startling global properties; for example, they can shrink curves and surfaces to round points \cite{Gage1986, Gage1984, Grayson1987}. These features made these flows powerful tools for proving topological theorems which ultimately led, via their generalizations by Hamilton \cite{Hamilton}, to Perelman's celebrated proof \cite{Perelman} of the Poincar\'e conjecture. The global features of curvature-driven flows are mostly related to the monotonic change of quantities, such as the entropy associated with Gaussian curvature \cite{Chow1}, other functionals, such as the Huisken functional \cite{Huisken1} in case of the Mean Curvature flow, the number of critical points of the curvature (used in the Curvature Scale Space model for image processing \cite{Mokhtarian,Mokhtarian1}), or the number of spatial critical points (with respect to a chosen reference point) \cite{Domokos2015,Grayson1987}, which are closely related to the geometry of the caustic \cite{Giblin1,Giblin2}.

Beyond offering powerful tools to prove mathematical statements, curvature-driven flows also have broad physical applications ranging from surface growth \cite{KPZ,Ishiwata2017} through image processing \cite{Koenderink,Mumford} and modeling cell migration and chemotaxis \cite{MacDonald2016} to mathematical models of abrasion \cite{Bloore1977,Firey1974}. Our paper is primarily motivated by the latter applications, where curvature-driven flows represent the fundamental model for the abrasion of pebbles under impacts of large particles. In particular, equation \ref{eq:curveshortening} has been proposed by Firey \cite{Firey1974} to model the extreme case of shape evolution of pebbles under collisions with infinitely large abraders. The other extreme scenario is abrasion by infinitely small particles, which is modeled by the so-called \emph{Eikonal equation}
\begin{equation}
\label{eq:eikonal}
v \equiv 1,
\end{equation}
also called the parallel map, arising in the study of wave fronts with constant speed that satisfy Huygens's principle.

Bloore \cite{Bloore1977} showed that the governing PDE of the general abrasion  model, including abraders of arbitrary size, is a linear combination of \ref{eq:curveshortening} and \ref{eq:eikonal}, resulting in the equation
\begin{equation}
\label{eq:bloore1}
v=1+C\kappa,
\end{equation}
which is often referred to as the (planar) Bloore flow \cite{DomokosG2012}. Moreover, 
the curve-shortening flow \ref{eq:curveshortening} has been generalized in a purely mathematical context, dominantly through the work of Andrews \cite{Andrews2002b, Andrews2002a, Andrews2002}, to the study of flows of the form
\begin{equation}
\label{eq:andrews1}
v=C\kappa^{\alpha}
\end{equation}
that are driven by a given power $\alpha$ of the curvature. 

We are concerned here with a two-parameter family that encompasses all of the above flows, which we write in a conveniently rescaled form as
\begin{equation}
\label{eq:abflow}
v = c + \kappa^\alpha.
\end{equation}
Here $c$ and $\alpha$ determine the relative strengths of the constant and curvature terms. Note that the constant $c$ in \ref{eq:abflow}, which we henceforth refer to as the Andrews-Bloore flow,  should not be confused with the constant $C$ in  \ref{eq:curveshortening}, \ref{eq:bloore1} and \ref{eq:andrews1}, which can be scaled to 1 by linearly scaling the curve $\Gamma$. 

Indeed, with this rescaling, we recover from \ref{eq:abflow}: the curve-shortening flow \ref{eq:curveshortening} for $c=0$ and $\alpha=1$; the Eikonal flow \ref{eq:eikonal} for $c=\alpha=0$; the Bloore flow \ref{eq:bloore1} for $c = \alpha = 1$; and the Andrews flow for $c=0$. We already also mention the special case $c=0$ and $\alpha=\frac{1}{3}$ of \ref{eq:abflow}, known as the \emph{affine curve-shortening flow} \cite{Alvarez}, which we will consider in section \ref{sec:evolellipse}. 

\subsection{Classical geophysical shape descriptors}

 In geophysical applications there has been particularly keen interest in identifying shape descriptors which, on the one hand, can be reliably measured and, on the other hand, evolve in a monotonic manner under these flows. The latter property makes the time evolution invertible, enabling geophysicists to deduce the provenance of pebbles based on observation of their current shape \cite{Szaboetal}. 
 
 One well-known and broadly applied shape descriptor is \emph{roundness}, commonly identified by the isoperimetric quotient $Q=4\pi A/L^2$; here $A$ is the enclosed area and $L$ the length of its perimeter, that is, the total arclength of the curve $\Gamma$. The monotonicity of $Q(t)$ under \ref{eq:curveshortening} was proven by  Gage \cite{Gage1984}; Andrews \cite{Andrews1994} generalized this result to \ref{eq:andrews1} and showed that the parameter value $\alpha=1/3$ (defining the affine shortening flow) separates flows with monotonically increasing and monotonically decreasing evolutions for $Q(t)$. We will illustrate these results in section \ref{sec:andrewsalpha}, where we not only illustrate the monotonicity of $Q(t)$ but also that, depending on $\alpha$, it approaches nontrivial limits predicted by Andrews \cite{Andrews2002b}.

Roundness $Q$ is a classical geophysical shape descriptor which is relatively easily measured \cite{Szaboetal} and, as described above, displays monotonic evolution under \ref{eq:curveshortening}. Still, it has a major handicap: under the Eikonal flow \ref{eq:eikonal}, which also serves as a relevant abrasion model in its own right, $Q(t)$
is \emph{decreasing} monotonically \cite{DomokosLangi2019}, so the direction of its evolution is \emph{opposite} to the evolution under purely curvature-driven abrasion. Hence, if both types of abrasion occur, the evolution of $Q(t)$ is unpredictable and unreliable.

\subsection{Critical points and their evolution}
\label{sec:critpoints}

To overcome the above difficulty, the (integer) number of critical points has been introduced \cite{Domokos2019, pebbles,VarkonyiDomokos2006} as a new type of shape descriptor. In this paper we will perform computations related to four alternative versions of critical points. In order to introduce them we first define some concepts and notation.

\begin{definition}  
\label{def_notation}
We consider the following objects related to planar curves: \\[-4mm]
\begin{itemize}
\item[(a)] 
Throughout, we assume that any curve $\Gamma$ is closed and simple, so that it encloses a topological disk of area $A$; the center of mass (or centroid) of this shape with homogeneous mass distribution is denoted by $C$.\\[-4mm]
\item[(b)] 
For $0 < \alpha$, the flow \ref{eq:abflow} reduces the area $A(t)$  enclosed by the evolving curve $\Gamma(t)$ in such a manner that $A(t)=0$ is reached for some finite time denoted $t = t_{\max}$.  We remark that for the curve-shortening flow \ref{eq:curveshortening} we have $t_{\max}=A(0)/2\pi$ \cite{Gage1986}; and for the Eikonal flow \ref{eq:eikonal} we have $t_{\max}=R_I$ \cite{DomokosLangi2019}, where $R_I$ is the radius of the largest inscribed circle of $\Gamma(0)$. \\[-4mm]
\item[(c)] 
The flow \ref{eq:curveshortening} shrinks the area $A(t)$ to a single point, called the \emph{vanishing point} or \emph{ultimate point}, which we denote by $U$. Whether or not $U$ exists for the flow \ref{eq:abflow} with $\alpha\neq1, c\neq 0$ is not clear \cite{Andrews2002b}. \\[-4mm]
\item[(d)] 
We further assume that any curve $\Gamma$ is star-like, which means that one can choose a point in its interior, called a reference point, with respect to which $\Gamma$ can be represented in polar coordinates; see \ref{sec:propGamma} and \ref{sec:polarmap} for details. Such a polar representation takes the from of a positive periodic radial distance function $r(\varphi)$ defined for $r(\varphi) \in [\-\pi,\pi]$; we stress that the distance function $r$ depends on the chosen reference point.
\item[(e)]  
The (signed) curvature of $\Gamma$, denoted by $\kappa$, can be represented similarly by a periodic function $\kappa(\varphi)$, which is computed from $r$ and its derivatives $r'$ and $r''$ with respect to $\varphi$; see \ref{sec:propGamma}. Throughout, we require that the curve $\Gamma$ is sufficiently smooth so that all required derivatives exist.
\end{itemize}
\end{definition}
We can now list the four types of critical points whose trajectories we will investigate during evolution of a given curve under the Andrews-Bloore flow \ref{eq:abflow}:
\smallskip
\begin{itemize}
\item[$n_{\mbox{\tiny $O$}}(t)$] 
denotes the number of critical points of $r(\varphi) = r_{\mbox{\tiny $O$}}(\varphi)$, characterized by $r' = dr/d\varphi=0$, when the reference point $O$ is chosen arbitrarily and remains fixed throughout the evolution;\\[-4mm]
\item[$n_{\mbox{\tiny $C$}}(t)$] 
denotes the number of critical points of $r(\varphi) = r_{\mbox{\tiny $C$}}(\varphi)$, characterized by $r' = dr/d\varphi=0$, when the reference point is the centroid $C(t)$ during the evolution, the location of which is generally not fixed under the flow; \\[-4mm]
\item[$n_{\mbox{\tiny $U$}}(t$)] 
denotes the number of critical points of $r(\varphi) = r_{\mbox{\tiny $U$}}(\varphi)$, characterized by $r' = dr/d\varphi=0$, when the reference point is the ultimate point $U$ throughout the evolution; \\[-4mm]
\item[$n_{\kappa}(t)$] 
denotes the number of critical points of $\kappa(\varphi)$ characterized by $d\kappa/d\varphi=0$; this number is independent of the chosen reference point for the periodic function $r$ used to compute the evolution.
\end{itemize}
\smallskip
In general, $n_{\mbox{\tiny $O$}}$ and $n_{\mbox{\tiny $U$}}$ have no direct physical meanings. On the other hand, $n_{\mbox{\tiny $C$}}$ is the number of static balance points of the homogeneous disk  bounded by the curve $\Gamma$, rolling along its perimeter on a frictionless, horizontal surface. 

In this paper we only treat the two-dimensional problem directly. However, we mention that the definitions are analogous in three dimension. In particular, static balance points of pebbles are easily determined either by hand experiments or by computer analysis of scanned image data. The most remarkable property of the different types of numbers of critical points above is that they appear to be decreasing monotonically (either in a deterministic or in a statistical sense) under all types of natural abrasion \cite{Domokos2019, DomokosLangi2019}. The main goal of our paper is to demonstrate this general property. This motivated the numerical algorithm presented here, designed for the purpose of tracking critical points during curve evolution. It allows us, on the one hand, to illustrate rigorously proven results regarding the evolution of critical points and, on the other hand, to perform computations in support of specific conjectures. We proceed by providing an overview of both kinds of results.

\subsubsection{Illustration of known results}
\label{sec:known}

The first rigorous result we consider is due to Grayson \cite{Grayson1987} who proved that under \ref{eq:curveshortening}, the number $n_{\mbox{\tiny $O$}}(t)$ is monotonically decreasing for any chosen, fixed reference point $O$. We will illustrate this phenomenon in \ref{sec:csf_trajectories}. Grayson's result was generalized in \cite{Domokos2015} to arbitrary curvature-driven flows $v=f(\kappa)$ under the condition that $df/d\kappa >0$ and, as one can observe, both the Bloore flow \ref{eq:bloore1} and the Andrews flow \ref{eq:andrews1} for $\alpha >0$ meet this criterion. 

Bloore proved \cite{Bloore1977} that the circle is an attractor for \ref{eq:bloore1} and, in the same paper, he investigated the limit as the curve $\Gamma$ approaches the circle, making an (albeit indirect) statement on the ultimate values of the numbers $n_{\mbox{\tiny $O$}}(t_{\max})$ and $n_{\mbox{\tiny $C$}}(t_{\max})$. As long as the limit shape is a round point, that is, the limit curve is a circle, the limit $n_{\mbox{\tiny $O$}}(t_{\max})=2$ appears to be obvious. Namely, except for the ultimate point $U$, any reference point $O$ fixed at $t=0$ will become external to the evolving curve for some $0 < t^{\star} < t_{\max}$. As the curve shrinks to a point, the ratio of its maximal diameter to its minimal distance to $O$ will approach zero; in this process $n_{\mbox{\tiny $O$}}$ will ultimately drop to $n_{\mbox{\tiny $O$}}=2$ and the two remaining trajectories of critical points will meet at the ultimate point $U$ at time $t=t_{\max}$ at 180 degrees. We will illustrate this property in \ref{sec:csf_trajectories} as well. Note that by data post-processing the trajectories of critical  can be determined also when $O$ is outside $\Gamma(t)$ and $r(\varphi)$ is no longer defined for all values $\varphi \ in [-\pi,\pi]$.

Bloore's stability theorem \cite{Bloore1977} claims that under \ref{eq:bloore1}, as $\Gamma$ approaches the circle while approaching the ultimate point $U$, Fourier terms vanish in reverse order, that is, the lowest-order terms decay at the slowest rate. This \emph{suggests}, at least in the absence of any symmetry (the generic case), that $n_{\mbox{\tiny $C$}}(t_{\max})=4$.  Moreover, we expect that the four trajectories of $n_{\mbox{\tiny $C$}}(t)$ that remain, meet at $U$ at right angles, that is, with 90 degree spacing. We remark that since $4 \leq n_{\mbox{\tiny $C$}}$ \cite{Ruina}, this statement also means that for generic curves the absolute minimum is achieved. This phenomenon is illustrated and confirmed for some examples evolving under the curve-shortening flow in \ref{sec:csf_trajectories}.

We will also illustrate interesting phenomena related to  the case with $\alpha<1$. The limiting shape of a simple closed planar curve under the evolution of the Andrews-Bloore flow~\ref{eq:abflow} for $c=0$ (that is, without a constant term) has been examined for different values of the exponent $\alpha$ of the curvature. For the case $\alpha = 1$ of curve-shortening flow, the curve converges to a point and the limiting shape is a circle; this was proven for convex curves by Gage~\cite{Gage1984,Gage1986} and non-convex curves by Grayson~\cite{Grayson1987}, implying also that the isoperimetric quotient $Q$ converges to 1.  Andrews examined how this limit depends on the exponent $0 < \alpha$ \cite{Andrews1996,Andrews1998,Andrews2002a}. He showed that for $\frac{1}{3} < \alpha $ any smooth convex curve converges to a point, the limiting curve is a circle and, hence, the isoperimetric quotient converges to 1. For $\alpha = \frac{1}{3}$ (the case of the affine shortening flow), the ellipses are homothetically contracting solutions, maintaining their axis ratios and isoperimetric quotients during the evolution. For $0 < \alpha \leq \frac{1}{3}$ the circle is not a limit any more and for any generic curve the isoperimetric quotient $Q$ approaches zero as $t \to t_{\max}$ \cite{Andrews2002a}. We will present illustrations of this dependence on $\alpha$ in \ref{sec:andrewsalpha}.

In subsequent work, Andrews \cite{Andrews2002} provided further insight into \ref{eq:abflow} with $c=0$. In particular, he showed that there are additional homothetic embedded curves with $D_n$-symmetry, which appear one-by-one as $\alpha$ is decreased down to zero. We illustrate and explore this property of the Andrews flow by presenting evolutions of $D_n$-symmetric curves in \ref{sec:highersymm}, where the exact statements and conditions on $\alpha$ can be found. Moreover, we demonstrate in \ref{sec:Cncurves} the result from \cite{Andrews2002b} that curves with $C_n$-symmetry converge to the homothetic solutions with $D_n$-symmetry in the $\alpha$-range where the latter exist.

\subsubsection{Conjectures supported by numerical evidence}
\label{sec:conj}

From the known results above, it appears as though for the case of the curve-shortening flow \ref{eq:curveshortening} the number of critical points of $r(\varphi)$ of a curve \emph{either} decreases monotonically when the reference point is a fixed chosen point \emph{or}, when the reference point does not involve an arbitrary choice but is taken to be the centroid, monotonicity cannot be proven. However, it is not quite the case that these two situations are mutually exclusive: for $n_{\mbox{\tiny $U$}}(t)$ we have an evolution that satisfies both criteria simultaneously. Since  $U$ is fixed, the results from \cite{Grayson1987, DomokosG2012} apply and thus $n_{\mbox{\tiny $U$}}(t)$ will behave like $n_{\mbox{\tiny $O$}}(t)$ and decrease monotonically. On the other hand, as far as the limiting behaviour at $t=t_{\max}$ is concerned, we suspect that it is analogous to that of $n_{\mbox{\tiny $C$}}(t)$. We formulate this statement as follows:
\begin{conjecture}
\label{conj_U}
Let the generic curve $\Gamma$ evolve under the curve shortening flow
\ref{eq:curveshortening}. Then we have $n_{\mbox{\tiny $U$}}(t_{\max})=4$ and the four remaining trajectories meet at right angles at $U = C(t_{\max})$.
\end{conjecture}
We will provide numerical evidence for this conjecture in \ref{sec:csf_trajectories}. Despite the obvious and conjectured good features of $n_{\mbox{\tiny $U$}}(t)$, from the practical point of view, it does not appear to be an attractive candidate to track erosion, because one can only compute $n_{\mbox{\tiny $U$}}(t_0)$ at a time $0 \leq t_0 \leq t_{\max}$ of interest after having computed the entire evolution up to $t_{\max}$, as required to identify the reference point $U$ \cite{Bryant1995}.  

We finally mention the last type of critical point, which does not have this disadvantage --- the number $n_{\kappa}(t)$ of the critical points of the (signed) curvature $\kappa$. Indeed, this number can be computed for any given curve; moreover, $n_{\kappa}$ is defined in an intrinsic manner and does not depend on the choice of the reference point. It has not been proven whether or not  $n_{\kappa}(t)$ decreases monotonically under evolution, but we suggest that it does:
\begin{conjecture}
\label{conj_kappa}
Let the generic curve $\Gamma$ evolve under the curve shortening flow \ref{eq:curveshortening}. Then $n_{\kappa}(t)$ decreases monotonically, and we have $n_{\kappa}(t_{\max})=4$ and the four remaining trajectories meet at right angles at $U$.
\end{conjecture}
We remark that due to the four-vertex theorem \cite{Kneser} we have $ n_{\kappa} \geq 4$, so our conjecture claims that the absolute minimum is reached, namely for all generic curves. In \ref{sec:csf_trajectories} we will also provide numerical evidence for this conjecture. 

Although $n_{\kappa}(t)$ appears to have attractive features, its application in geomorphology is limited. First of all, determining $n_{\kappa}(t_0)$ at any given time $0 \leq t_0 \leq t_{\max}$ of interest is not possible by hand experiments. Moreover, from the computational point of view, it is more challenging to determine $n_{\kappa}$ from a scanned image data than determining either $n_{\mbox{\tiny $O$}}$ or $n_{\mbox{\tiny $C$}}$; this is the case because it requires third-order derivatives of the distance function $r(\varphi)$, rather than second derivatives.

\subsection{Essence of the algorithm}

We present in this paper a numerical method for computing the evolution under the flow \ref{eq:abflow} of a given curve in the plane --- designed specifically to reliably track (the trajectories of) critical points along the evolution, all the way up to (very close to) the final point. We require that the initial curve admits polar coordinates from a point in its interior. Such curves are known as star-like, and $r(\varphi)$ is given as the distance from the reference point of the unique intersection point of the curve with the ray of angle $\varphi \in [-\pi,\pi]$. Note further that star-like curves are a natural generalisation of convex curves, and that the flow \ref{eq:abflow} generally preserves the star-like nature of the initial curve; in fact, it is typically the case that a star-like curve becomes convex under evolution.

The algorithm's capability of tracking critical points of the radial distance function $r(\varphi)$ is achieved by representing $r(\varphi)$ by a piecewise polynomial (over a specified mesh). This is natural in light of our wish to study the functions $n_{\mbox{\tiny $O$}}(t)$, $n_{\mbox{\tiny $C$}}(t)$ and $n_{\mbox{\tiny $U$}}(t)$ of critical points, as determined by the chosen fixed or moving reference point, as well as the number of critical points $n_{\kappa}(t)$. The general setup of our method is quite flexible and allows for different types of meshes and remeshing strategies. The key property of our approach is that every curve during the evolution is represented by a sufficiently smooth function that can be evaluated for any angle, that is, anywhere along the curve and not just at the mesh points. Since our representation is polynomial, the derivatives of the radial distance function $r$, the normal and the curvature $\kappa$ at any point of the curve, as well as certain observables, can be computed symbolically, that is, do not need to be approximated by using numerical differentiation. 

Indeed, we will demonstrate that, with appropriately chosen accuracy settings, our algorithm is able to track reliably and with manageable computational effort the critical points of the radial distance function $r$ and of the curvature $\kappa$ during evolutions of curves under \ref{eq:abflow}, including for choices of the parameters $c$ and $\alpha$ beyond the standard case $c=0$ and $\alpha=1$ of the curve-shortening flow. It is this new capability that distinguishes our formulation from existing methods for computing curvature-dependent flows. Reaction-diffusion techniques \cite{Barrett2017}, level-set method \cite{Osher1988,Osher2002,Sethian1999}, and front-tracking methods based on finite differences \cite{Dziuk1999,Mackenzie2019} and related finite element methods \cite{Barrett2011,Dziuk1999,MacDonald2016,Mikula1999} have been used for this purpose. These computational approaches are highly efficient in several situations. Generally, the focus has been on some particular property of the problem or application, such as the curve reaching a singularity and breaking up into a set of disconnected curves in the case of level set methods. Moreover, front-tracking methods, in particular, which may also feature different types of mesh adaptation  \cite{Barrett2011,Mackenzie2019}, can be used to visualize the evolution of curves (that need not be star-like or convex) under the curve-shortening flow \ref{eq:curveshortening}; see, for example, the illustrations in \cite{Daskalopoulos2020} of results from \cite{Osher1988,Malladi1997}. On the other hand, the evolution of shape descriptors, especially trajectories of critical points and their numbers $n_{\mbox{\tiny $O$}}(t_{\max})$ and $n_{\mbox{\tiny $C$}}(t_{\max})$, has not been the focus of computations before. At the early stages of this project, we found that the existing methods are not well suited to capturing the fine details of the curve required for this task: even with very fine meshes and/or extremely small time steps it is difficult to achieve and maintain the required accuracy as the final point is approached. The approach presented here is complementary to available methods and has been inspired by collocation methods for boundary value problems of differential equations \cite{doedel2007,iserles1996a}. The underlying idea is to work with a space of smooth approximating functions \cite{rivlin,trefethenbook}, generally piecewise polynomials, to represent the global objects under consideration such as periodic orbits of differential equations or curves or fronts evolving under a geometric flow. Another important aspect of our method is that the planar curve is represented by a function of a single variable, namely the periodic radial distance function $r$, whose critical points are therefore readily available. While this requires that any curve be star-like, this is actually quite a natural assumption in the context of abrasion problems and limiting shapes under curvature flows. A method that is quite similar in spirit is that for front propagation in the Navier-Stokes equations in \cite{Popinet1999}, where the front is approximated with a parametric representation based on polynomials over a list of points; the front is then propagated by interpolating the velocity field over these points.

\subsection{Structure of the paper}

Section \ref{sec:propGamma} first introduces the polar parameterization of a star-like curve and then details how all relevant quantities can be computed from its periodic radial distance function. We then discuss the beneficial properties of the central parameterization (with the centroid as the reference point) in \ref{sec:Cparam}, how symmetries of a curve are reflected in its polar parameterization in \ref{sec:symmGamma}, and properties of trajectories of critical points in \ref{sec:evolution}. Our computational method is then presented in \ref{sec:alg}, where we introduce the space of discretized curves in \ref{sec:discr_curve} and then formulate the general algorithm in \ref{sec:discr_int}; time stepping, the choice of mesh, its propagation and remeshing are discussed in \ref{sec:timestep} through \ref{sec:remesh}, respectively. Section \ref{sec:accuracy} then provides a number of benchmark test that demonstrate the effects of different accuracy settings: the mesh size in \ref{sec:initialproj} and \ref{sec:evolmesh}, and the type of mesh and the remeshing strategy in \ref{sec:evolremesh}; we then show in \ref{sec:rotatedsymm} that trajectories of critical points can be computed reliably, irrespective of the position of the curve; the final \ref{sec:evolellipse} shows that an initial ellipse does indeed not change its shape (within the numerical accuracy of the computation) under the affine shortening flow, as required by a result of Andrews \cite{Andrews1996}. In \ref{sec:examples} we then provide two case studies that demonstrate how our algorithm can be used as a tool for the experimental investigation of the evolution under \ref{eq:abflow} of different types of curves; in \ref{sec:csf_trajectories} we provide numerical evidence in support of \ref{conj_U} and \ref{conj_kappa}. In the final \ref{sec:andrewsalpha} we illustrate results of Andrews regarding the shape evolution for $c=0$ and $0 < \alpha \leq 1$ and investigate the limit shapes of curves with $D_n$-symmetry and with $C_n$-symmetry. In \ref{sec:conclusions} we draw some conclusions and point to future work. Finally, \ref{sec:polarmap} provides some additional background on the existence of a polar parameterization.

\section{Polar parameterization and properties of $\Gamma$}
\label{sec:propGamma}

Throughout, we work within the set of star-like planar curves; the shape bounded by such a curve is referred to as a star domain or radially convex set. Any star-like curve $\Gamma$ is closed and simple and further characterized by the property that we can find a point $O$, referred to as a \textit{reference point}, such that all rays from $O$ intersect $\Gamma$ transversely. Since this property is open, there is an open set of possible reference points in the interior of $\Gamma$. Note any convex curve is star-like with the additional and defining property that any point in its interior can be chosen as a reference point. The latter is not the case for non-convex star-like curves, which necessarily have points of inflection, that is, points where the curvature changes sign; see \ref{sec:polarmap} for further details.

For any star-like planar curve $\Gamma$ and chosen reference point $O$ the polar parameterization around $O$ takes the form 
\begin{equation}
\label{eq:gammapol}
\begin{array}{rcl}
\gamma: [-\pi,\pi]  & \to & \bbR^2 \cong \bbC \\
\varphi \ \ \ & \mapsto & O + r(\varphi) e^{i\varphi} 
\end{array}
\end{equation}
where 
\begin{equation}
\label{eq:def_r}
\begin{array}{rcl}
r: [-\pi,\pi] & \to & \bbR^+\\ 
    \varphi \ & \mapsto & r(\varphi) \quad {\rm with} \quad 
O + r(\varphi) e^{i \varphi} \in \Gamma .
\end{array}
\end{equation}
Note that we identify $\bbR^2$ with $\bbC$ to allow for compact notation involving multiplication in $\bbC$ where convenient. Hence, $r(\varphi)$ is the radial distance from $O$ in the direction $\varphi$, where the angle $\varphi \in [-\pi, \pi]$ is measured in radians in the usual way (from the positive $x$-direction from the reference point $O$ and in the mathematically positive direction). In other words, $\varphi$ ranges over the closed fundamental interval $[-\pi,\pi]$ of the covering space $\bbR$ of the circle $2\pi\bbS^1 = \bbR/2\pi\bbZ$; we inlude both $-\pi$ and $\pi$ in the domain to stress the periodicity of $r$ given by $r(-\pi)=r(\pi)$. By the very definition of $\gamma$ the curve $\Gamma$ is given by
\begin{equation}
\label{eq:Gammapol}
\Gamma  = \{\gamma(\varphi) \ | \ \varphi \in  [-\pi,\pi] \} .
\end{equation}
Clearly, the polar parameterization and the functions $r=r_{\mbox{\tiny $O$}}$ and $\gamma=\gamma_o$ depend on the choice of reference point $O$; in what follows, we denote the dependence on $O$ in the notation only where this is of specific importance. 

We remark for future reference that an arclength parameterization $\tilde{\gamma}(s)$ of $\Gamma$ can be obtained when needed from the polar parameterization \ref{eq:gammapol} as
\begin{equation}
\label{eq:polartoarc}
\tilde{\gamma}(s) := \gamma(\varphi) \quad {\rm where} \  \varphi = \varphi(s) \ 
{\rm is \ given \ by} \quad  
s = \int_{-\pi}^{\varphi(s)} \| \gamma'(\theta) \| \, d\theta .
\end{equation}
Here $0 \leq s \leq L$ is the arclength parameter and $L$ is the total arclength of $\Gamma$.

The periodic radial distance function $r$ from \ref{eq:def_r} is well defined for any star-like curve and contains all information on the curve $\Gamma$. In particular, we can express the parameterizations of the normal $\n$ and the curvature $\kappa$, both needed to evaluate the flow \ref{eq:abflow}, as periodic functions of $\varphi$, in terms of $r$ and its first and second derivatives $r'$ and $r''$ with respect to $\varphi$ as follows. As mentioned in \ref{def_notation}(e), we assume throughout that $\Gamma$ is smooth enough to allow us to take any derivatives needed. With
\begin{equation}
\gamma'(\varphi) = (r'(\varphi) + i r(\varphi) ) e^{i\varphi} \quad {\rm and} \quad \
\end{equation}
\begin{equation}
\gamma''(\varphi) = (r''(\varphi) + 2 i r' (\varphi) - r (\varphi)) e^{i\varphi} 
\end{equation}
the unit normal vector (perpendicular to the tangent and pointing inwards) is given by 
\begin{equation}
\label{eq:def_n}
\n(\varphi) = \frac{i \gamma'(\varphi) }{\| \gamma'(\varphi) \|}  
\end{equation}
and the signed curvature by
\begin{equation}
\label{eq:def_kappa}
\kappa(\varphi) = \frac{{\rm det}[ \gamma'(\varphi) \,\, \gamma''(\varphi)]}
{\| \gamma'(\varphi) \|^3} .
\end{equation}

It is insightful to write the functions $\n$ and $\kappa$ in polar coordinates around $O$; here we do not show for notational convenience the dependence on $\varphi$ of the periodic functions introduced above. First of all,
\begin{equation}
\label{eq:polparphi}
\gamma'
=  \| \gamma' \| e^{i(\varphi + \beta)} 
\end{equation}
where the angle $\beta = \beta(\varphi)$ (with respect to the positive $x$-axis as usual) satisfies
\begin{equation}
\label{eq:propalpha}
\cos(\beta)
= \frac{r'}{\|\gamma' \|}
= \frac{r'}{\sqrt{r^2 + (r')^2}}
\quad {\rm and} \quad
\sin(\beta)
= \frac{r}{\|\gamma' \|}
= \frac{r}{\sqrt{ r^2 + (r')^2}}
\end{equation}
Similarly,
\begin{equation}
\gamma'' 
=  \| \gamma'' \| e^{i(\varphi + \eta)} 
\end{equation}
where the angle $\eta = \eta(\varphi)$ 
satisfies
\begin{equation}
\label{eq:propbeta}
\cos(\eta)
= \frac{r'' - r}{\| \gamma'' \|}
= \frac{r'' - r}{\sqrt{(r'' - r)^2  + 4 (r')^2}}
\quad {\rm and} \quad
\sin(\eta)
= \frac{2 r'}{\| \gamma'' \|}
= \frac{2 r'}{\sqrt{(r'' - r)^2  + 4 (r')^2}}.
\end{equation}
This gives
\begin{equation}
\label{eq:def_n_angle}
\n= \frac{i \gamma' }{\| \gamma' \|}  
= e^{i(\varphi + \beta + \frac{\pi}{2})}  
\end{equation}
and
\begin{equation}
\label{eq:def_k_r}
\kappa  
 =  \frac{\| \gamma'' \|}{\| \gamma' \|^2} 
\sin(\eta - \beta)
= \frac{r(r -r'') + 2(r')^2}{(r^2 + (r')^2)^\frac{3}{2}} .
\end{equation}
It follows from \ref{eq:def_n_angle} that the normal $\n$ points at the reference point $O$ exactly at the critical points of $r$; namely, then $\n = -e^{i\varphi^*}$ and $\beta(\varphi^*) = \frac{\pi}{2}$, which is equivalent to $\cos(\beta(\varphi^*)) = 0$ and, hence, to $r'(\varphi^*) =0$ according to \ref{eq:propalpha}. 

Moreover, we can conclude from \ref{eq:def_k_r} that $\Gamma$ is convex, that is, does not have inflection points where $\kappa(\varphi^*) = 0$ for some $\varphi^*$, if and only if
\begin{equation}
\label{eq:condincl}
r^2 + 2(r')^2 >  r  r'', 
\end{equation}
which means that the graphs of the periodic functions $r^2 + 2(r')^2$ and $r  r''$ do not intersect. Since convexity is an intrinsic property of the curve $\Gamma$, 
condition \ref{eq:condincl} does not depend on the reference point; that is, it can be verified for the radial distance function of any reference point.  Note also that, $r'' < r$ implies \ref{eq:condincl} (but is not necessary); hence, if the graphs of $r$ and $r''$ do not intersect then $\Gamma$ is convex.

\subsection{The central parameterization}
\label{sec:Cparam}

As we discussed in \ref{sec:propGamma}, the radial distance function $r = r_{\mbox{\tiny $O$}}$ for any choice of reference point $O$ contains all of the information on $\Gamma$. So in this sense, any reference point is as good as any other reference point. However, in the context of shape evolution it is a natural choice to use the center of gravity or centroid $C$ of the (homogeneous) shape or region bounded by $\Gamma$ as the reference point. The centroid is intrinsic to the curve and can be computed from any polar parameterization as 
\begin{equation}
\label{eq:computeC}
C = \frac{1}{A} \int_{-\pi}^{\pi} \int_{0}^{r(\theta)} e^{i\theta} r^2\, dr d\theta
= \frac{1}{3A} \int_{-\pi}^{\pi}  r^3(\theta) e^{i\theta} d\theta.
\end{equation}
Here, $A$ is the area of the shape bounded by $\Gamma$, which is also intrinsic and given by 
\begin{equation}
\label{eq:computeA}
A = \int_{-\pi}^{\pi} \int_{0}^{r(\theta)} r\, dr d\theta
= \frac{1}{2} \int_{-\pi}^{\pi} r^2(\theta) d\theta.
\end{equation}
We refer to the specific polar parameterization with $O = C$ as the \textit{central parameterization} and to the periodic function $r$ as the \textit{central radial distance function}, which we denote $r_{\mbox{\tiny $C$}}$ where the context calls for it. The central parameterization has a number of advantages. \\[-3mm]

\begin{enumerate}
\item
The central parameterization encodes the relevant properties of $\Gamma$ in a particularly convenient way. Namely, the stationary points of the associated homogeneous shape are immediately available as the critical points of the central radial distance function $r_{\mbox{\tiny $C$}}$, as is the number of stationary points $n_{\mbox{\tiny $C$}}$; note that stationary points along negatively curved parts of $\Gamma$ need to be balanced `on a pin'; see \ref{sec:polarmap}. \\[-3mm]
\item
For the central parameterization any discrete symmetries of the curve $\Gamma$ are conveniently represented as symmetires of the central radial distance function $r_{\mbox{\tiny $C$}}$; see \ref{sec:symmGamma} for details. \\[-3mm]
\item
The reference point $C$ is inherent to the curve $\Gamma$. Hence, the central parameterization, if it exists, is a well-defined and definite choice of polar parameterization, which is an advantage from the algorithmic point of view. We remark that the centroid $C$ is generally not fixed under curvature flow, but can be tracked as $C(t)$ at negligible extra cost by evaluating the integral \ref{eq:computeC}; in particular, $C(t)$ converges to the ultimate points $U$ as $t \to t_{\max}$. \\[-3mm]
\item
Should information on the evolution of critical points of the radial distance function $r_{\mbox{\tiny $O$}}$ for any other reference point $O$ be required, it can be found from the central parameterization by post-processing; in other words, it is not necessary to recompute the evolution. \\[-3mm]
\end{enumerate}
For these reasons we will work in the implementation of our algorithm for the evolution of curves under \ref{eq:abflow} with the central parameterization whenever it exists, that is, when the centroid can be chosen as a reference point; see also \ref{sec:polarmap}. While this is not true in general, it is observed that for many choices of $c$ and $\alpha \geq 0$ the curvature flow \ref{eq:abflow} drives any star-like curve $\Gamma$ towards convexity. This has been proven by Grayson for $c=0$ and $\alpha=1$ \cite{Grayson1987} and by Andrews for certain intervals of $\alpha$ \cite{Andrews2002a}. In fact, for all our examples the centroid $C$ can be chosen as a reference point, so that the central parameterization is already available for the initial curve $\Gamma$.

\subsection{Symmetries of the curve $\Gamma$}
\label{sec:symmGamma}

We will also consider evolutions of curves that are invariant under a discrete symmetry group; typically, this will be the dihedral group $D_n$ of rotations and reflections of a regular $n$-gon in $\bbR^2$. So suppose that $\Gamma$ is invariant under the discrete group $G$, which means that $g \Gamma = \Gamma$ for all $g \in G$; here $g$ acts on $\Gamma$ as either a rotation or a reflection of $\bbR^2$. The question is whether these intrinsic symmetries of $\Gamma$ given by $G$ are reflected in the radial distance function $r_{\mbox{\tiny $O$}}$ of the polar parameterization with a given reference point $O$. This is only the case when $O$ is chosen to lie in the fixed point subspace of the respective group element or subgroup. More specifically, suppose $O \in {\rm Fix}(H)$ for some subgroup $H \leq G$, then $H$ induces an invariance of the scalar function $r_{\mbox{\tiny $O$}}$, namely
\begin{equation}
\label{eq:symmrotatephi}
r_{\mbox{\tiny $O$}}(\varphi) = r_{\mbox{\tiny $O$}}(\varphi + \frac{2\pi}{\ell})
\end{equation}
for a rotation over $\frac{2\pi}{\ell}$ about $O$ and 
\begin{equation}
\label{eq:symmreflectphi}
r_{\mbox{\tiny $O$}}(\theta + \varphi) = r_{\mbox{\tiny $O$}}(\theta - \varphi)
\end{equation}
for a reflection in a line through $O$ of angle $\theta$. Hence, if $O \in {\rm Fix}(G)$ then $r = r_{\mbox{\tiny $O$}}$ has all symmetries represented by $G$; indeed then $r_{\mbox{\tiny $O$}}'$ and $r_{\mbox{\tiny $O$}}''$ have the exact same symmetries of translation and reflection (in the exact same points). Note that the tangent $\n$ and the signed curvature $\kappa$ do not depend on the choice of reference point; hence, these periodic function of $\varphi$ are $G$-invariant in the same way (that is, under discrete translations and reflections) in any case.

The centroid $C$ of $\Gamma$ is fixed under all elements of $G$, that is, $C \in {\rm Fix}(G)$, so that the radial distance function $r_{\mbox{\tiny $C$}}$ has all symmetries of $G$. Moreover, unless $G = \bbZ_2$ is generated by a single reflection, $G$ contains a rotation, which must therefore be a rotation of $\Gamma$ around $C$. In particular, $C$ lies on the intersection of all lines of reflection. Note that the normal $\n(\theta)$ of a point $\gamma(\theta)$ on a line of reflection points to $C$; hence, generically, $r_{\mbox{\tiny $C$}}$ has a maximum or minimum, reflecting the fact that every point of $\Gamma$ on a line of reflection is a stationary point. It follows from \ref{eq:def_k_r} that the scalar function $\kappa$ has the same symmetry as $r_{\mbox{\tiny $C$}}$, which means that the stationary points on lines of reflection are also extrema of the curvature $\kappa$. The center $I$ of the largest inscribed circle must also be fixed under the action of the group. It follows that, if the group $G$ contains two reflections or a rotation, then $I = C$ and the extrema of the curvature $\kappa$ coincide with the stationary points.

Considering curves with discrete symmetry is interesting because the flow of \ref{eq:abflow} preserves symmetry; this follows from the fact that the normal and the curvature inherit the symmetry of $\Gamma$. Hence, for any $0 \leq t < t_{\max}$ the curve $\Gamma(t)$ is invariant under the symmetry group $G$ of $\Gamma$; we will see examples of this in the coming sections. In particular, any method to compute the evolution of a planar curve must preserve its symmetry. As we will see in \ref{sec:alg}, checking this property provides a good test of the performance of our implementation. Indeed, whether symmetry is preserved can be checked readily by considering the symmetry properties of the radial distance function $r_{\mbox{\tiny $C$}}$ under evolution. This is another reason why the central parameterization is a good choice of parameterization to work with algorithmically.

\subsection{Trajectories of critical points}
\label{sec:evolution}

The flow of equation \ref{eq:abflow} for given $c$ and $\alpha$ generates the (forward) evolution $\Gamma(t)$ of a given curve $\Gamma = \Gamma(0)$. We consider here only evolutions of star-like curves and compute the evolution $\Gamma(t)$ for $0 \leq t  \leq t_{\max}$ via the corresponding radial distance distance functions $r(\varphi;t)$ of its polar (and generally central) parameterization. For any $t$, the function $r(\varphi;t)$ is a periodic Morse function \cite{ArnoldC,Demazure,PostonStewart}, which means (assuming that $\Gamma(t)$ is sufficiently smooth) that the derivative $r'(\varphi;t) = \frac{d r}{d\varphi}(\varphi;t)$ vanishes at isolated points $\varphi^* =\varphi^*(t) \in [-\pi,\pi]$. Generically, these points correspond to an alternating sequence of maxima and minima separated by isolated roots of the second derivative $r''(\varphi;t) = \frac{d^2 r}{d\varphi^2}(\varphi;t)$; due to periodicity, there are generically equal numbers of maxima and minima, that is, an even number of extrema of $r(\varphi;t)$. It follows further from singularity theory \cite{ArnoldC,Demazure,GolSchaefer,PostonStewart} that the number of extrema changes generically at isolated points of the $t$-line. Namely, at such a point $t^*$ one has either \\[-3mm]
\begin{enumerate}
\item
In the absence of additional symmetry properties of $\Gamma(t)$: a fold bifurcation at $t^*$, which is a cubic singularity where $\frac{d r}{d \varphi}(\varphi^*;t^*) = \frac{d^2 r}{d \varphi^2}(\varphi^*;t^*) = 0$; the genericity condition is that $\frac{d^3 r}{d \varphi^3}(\varphi^*;t^*) \neq 0$. In this bifurcation a maximum and a minimum of $r$ disappear or appear as $t$ increases through $t^*$.  \\[-3mm]
\item
Across a line of reflection of $\Gamma(t)$: a pitchfork bifurcation at $t^*$ of an extremum of $r(\varphi^*;t^*)$, where $\frac{d r}{d \varphi}(\varphi^*;t^*) = \frac{d^3 r}{d \varphi^3}(\varphi^*;t^*) = 0$; here the genericity condition is that $\frac{d^5 r}{d \varphi^5}(\varphi^*;t^*) \neq 0$, and note that $r$ is an odd function with respect to $\varphi^* + \varphi$ due to reflectional symmetry. In this bifurcation a pair of minima/maxima off the reflection axis disappear or appear by meeting a maximum/minimum on the axis of reflection as $t$ increases through $t^*$.  \\[-3mm]
\item
With rotational symmetry of $\Gamma(t)$: a fold or pitchfork bifurcation at $t^*$ occurs simultaneously along the orbit of the generator of rotation, that is, for all $\frac{2\pi}{\ell} \varphi^*$ where $\ell$ is the order of the rotational subgroup of the symmetry group $G$ of $\Gamma(0)$.  \\[-3mm]
\end{enumerate}
The exact same statements hold also for the scalar signed curvature function $\kappa(\varphi;t)$ of $\Gamma(t)$, which is also a periodic Morse function determined from $r(\varphi;t)$ by \ref{eq:def_k_r}.  Since \ref{eq:abflow} preserves symmetry, it follow that the numbers $n_{\mbox{\tiny $O$}}(t)$, $n_{\mbox{\tiny $C$}}(t)$, $n_{\mbox{\tiny $U$}}(t)$ and $n_\kappa(t)$ are generically constant and even, and change at isolated points $t^*$ by $2 \ell$ where $\ell$ is the order of the rotation subgroup of the symmetry group of $\Gamma=\Gamma(0)$ (which may be empty, in which case $\ell = 1$). 

The challenge, from both the theoretical and the algorithmic point of view, is to accurately compute the trajectories of the different types of critical points during the evolution of a given curve $\Gamma$. The numerical method we introduce below is specifically designed for this task. In particular, it allows us to illustrate, confirm or make statements about the properties of the numbers $n_{\mbox{\tiny $O$}}(t)$, $n_{\mbox{\tiny $C$}}(t)$, $n_{\mbox{\tiny $U$}}(t)$ and $n_\kappa(t)$ as $t \to t_{\max}$.

\section{Computing the evolution of a curve}
\label{sec:alg}

For any star-like curve $\Gamma$ equation \ref{eq:abflow} assigns at every point $\gamma(\varphi)$ the infinitesimal velocity $v(\varphi) = (c + \kappa^\alpha (\varphi)) \n(\varphi)$. Here $\gamma$ is the polar parameterization of $\Gamma$ with respect to a suitable chosen reference point, with is the centroid $C$ in our implementation; hence, we do not indicate the dependence on the reference point of the different periodic functions in the notation from now on. Computing the evolution $\Gamma(t)$ under the flow given by \ref{eq:abflow} is therefore an initial-value problem with initial condition $\Gamma$. Solving this initial-value problem numerically requires a discretization of time by means of an integration step, as well as a discretization of space, that is, of the curve that is being evolved. 

Our goal is to compute, a sequence of curves $(\Gamma_{i})_{0 \leq i \leq I}$ up to some $I \in \bbN$ with central parameterizations $\gamma_{i}$ that provide an accurate approximation of the sequence $(\Gamma(t_i))$ with $t_i = \sum_{l=1}^{i} h_l$. Here, the $h_i$ with $\sum_{i=1}^{I} h_i < t_{\max}$ are the time steps of the integration during such a computation, which need to be chosen and adapted in a suitable way.

\subsection{The space of discretized curves}
\label{sec:discr_curve}

We first deal with the issue that curves in $\bbR^2$ are infinite-dimensional objects. Hence, the task is to discretize any star-like curve $\Gamma$ in an appropriate way. To achieve this we consider a finite-dimensional space $\mathscr{D}$ of discretized curves that are defined over a mesh of $N$ points. There are several choices one could make for the mesh and the space of functions that specify $\mathscr{D}$. We work here with a mesh 
\begin{equation}
\label{eq:def_gammamesh}
M := \{m_j\}_{0 \leq j \leq N} \subset \Gamma \subset \bbR^2
\end{equation}
of mesh points $m_j$ on the given curve $\Gamma$. Since $\Gamma$ has a given polar parameterization $\gamma$, the mesh $M$ can be pulled back to a mesh
\begin{equation}
\label{eq:def_phimesh}
M^\varphi := \{\varphi_j\}_{0 \leq j \leq N} \subset [-\pi, \pi]
\quad {\rm with} \quad \gamma(\varphi_j) = m_j
\end{equation}
of the polar angle $\varphi$. Reversely, for a polar-angle mesh $M_\varphi$ the parameterization $\gamma$ induces a unique mesh of points on $\Gamma$, and we express this relationship between the two dual meshes as
\begin{equation}
\label{eq:def_mapmesh}
M^\varphi = \gamma^{-1}(M) \quad {\rm and} \quad M := \gamma(M^\varphi)
\end{equation}
for notational convenience. 

We now consider piecewise polynomial radial distance functions through the points of the polar-angle mesh $M^\varphi$. Again, there are different options for choosing such polynomials. We define and compute them here as interpolating polynomials of degree $p$ on each mesh interval $[\varphi_i, \varphi_{i+1}]$, whose first $p$ derivatives agree at the mesh points of $M^\varphi$. We refer to this class of periodic piecewise polynomial functions as $\mathscr{P}_N^p$; they are also known as periodic splines of degree $p$ \cite{ANW67}.  This choice defines the discretization space 
\begin{equation}
\label{eq:def_subspace}
\mathscr{D} = \{ \mbox{\rm curves} \ \Gamma  \ 
\mbox{\rm having a polar parametrization \eqref{eq:gammapol} with} \ 
r \in \mathscr{P}_{N_{\max}}^p \}. 
\end{equation}
Here $N_{\max}$ and $p$ are taken to be fixed; we find that degree $p = 5$ polynomials are a suitable choice, in practice, in terms of the balance between providing sufficient smoothness versus the cost of computing the polynomial representation. Similarly, there is a trade-off between run-time and data size versus accuracy when it comes to the choice of the mesh size. Note that we require that the mesh size $N$ remains bounded by some apriori bound $N_{\max}$ solely to ensure that $\mathscr{D}$ is of finite dimension; hence, our setup allows for the adaptation of the mesh and its size during the computation of an evolution. We find that meshes of $N$ between 50 and 100 generally suffice for the curves we consider; see section \ref{sec:accuracy} for further details.

It is important to note that the radial distance function $r(\varphi)$ of any curve $\Gamma \in \mathscr{D}$ is defined for all $\varphi \in [-\pi, \pi]$. In other words, we are still working on a space of periodic functions and not just on a set of mesh points. In particular, the derivatives $r'$ and $r''$, as well as the normal $\n$ and curvature $\kappa$ are also piecewise polynomial functions, which can be computed symbolically from $r \in \mathscr{P}_{N_{\max}}^p$. Hence, all of these periodic functions can be evaluated exactly at any point $\varphi \in [-\pi, \pi]$. 

Any curve star-like curve $\Gamma$ is restricted to the discretization space $\mathscr{D}$ in a unique way by the choice of the mesh $M$ (and the reference point); we refer to this restriction or projection to the curve with the corresponding radial distance function $r \in \mathscr{P}_{N_{\max}}^p$ as $\Pi_M(\Gamma)$ for notational convenience. Note that the piecewise polynomial radial distance function $r$ of the projected curve $\Pi_M(\Gamma)\in \mathscr{D}$ can be computed readily for any chosen mesh and reference point; in what follows, the reference point is always the centroid C unless otherwise stated.

\subsection{General formulation of the curve evolution algorithm}
\label{sec:discr_int}

For a curve $\Gamma \in \mathscr{D}$ with a central parameterization $\gamma$ with piecewise polynomial radial distance function $r \in \mathscr{P}_{N_{\max}}^p$ for some given mesh $M$ we perform the integration with an Euler step $E_h$ of stepsize $h$ for the flow of \ref{eq:abflow} given by
\begin{equation}
\label{eq:Euler}
E_h(\gamma)  = \gamma + h v, 
\end{equation}
where
\begin{equation}
\label{eq:Eulervarphi}
E_h(\gamma) (\varphi) = C + e^{i\varphi} r(\varphi) + h (c + k ^\alpha (\varphi)) \n(\varphi)).
\end{equation}
The resulting, $h$-evolved periodic orbit $\Gamma_h = E_h(\Gamma)$ is, in general, not in the discretization space $\mathscr{D}$, that is, it does not have a polar parameterization with a piecewise polynomial radial distance function. This problem is overcome by constructing the new and unique radial distance function $r_h \in \mathscr{P}_{N_{\max}}^p$ and associated parameterization $\gamma_h(\varphi)$ with respect to a prespecified mesh $\widehat{M} \subset E_h(\Gamma)$ (and given reference point). Formally, this is the projection $\Pi_{\widehat{M}}$ from $\mathscr{C}$ back to $\mathscr{D}$, which can be expressed as 
\begin{equation}
\label{eq:EulerpojectD}
\gamma_h(\varphi) = \Pi_{\widehat{M}} (E_h(\gamma)) (\varphi)  = C + e^{i\varphi} r_h(\varphi).
\end{equation}
Therefore, the integration step on the discretization space $\mathscr{D}$ of a curve $\Gamma \in \mathscr{D}$ is defined as an Euler step followed by projection back onto $\mathscr{D}$, that is, by  $\Pi_{\widehat{M}} (E_h(\Gamma)) \subset \mathscr{D}$; note that this requires one to specify the mesh $\widehat{M}$ used in the restriction back to $\mathscr{D}$.

The evolution algorithm to approximate the evolution $\Gamma(t)$ of a given star-like curve $\Gamma$ under the flow of \ref{eq:abflow} now consists of computing a sequence $(\Gamma_{i})_{0 \leq i \leq I} \subset \mathscr{D}$ as
\begin{equation}
\label{eq:trajD}
\begin{array}{rcl}
\Gamma_0 & = & \Pi_{M_0}(\Gamma) \\ 
    \Gamma_{i}  & = & (\Pi_{M_i} \circ E_{h_i}) (\Gamma_{i-1}) .
\end{array}
\end{equation}
By construction, the parameterizations 
\begin{equation}
\label{eq:trajDgamma}
\gamma_i = (\Pi_{M_i} \circ E_{h_i}) (\gamma^{i-1}) 
\end{equation}
have central radial distance functions $r_i$ that are piecewise interpolating polynomials defined by the meshes $M_i$. Where it is appropriate to stress that the sequence $(\Gamma_{i})$ represents a time evolution, we also refer to $\Gamma_{i}$ as $(\Gamma_{t_i})$, where $t_i = \sum_{l=1}^{i} h_l$ is the integration time up to step $i$ of the computation.

Formulating the evolution algorithm in general terms first is convenient because it encompasses different strategies for updating the mesh, that is, the discretization of the curve that is being evolved. In the remainder of this section, we present our numerical implementation of the general evolution algorithm given by \ref{eq:trajD}. To this end, we will specify: (1) the stepsize sequence $(h_i)$ with suitable stopping criteria for the computation, and  (2) a strategy for choosing the mesh $(M_i)$ at step $i$ of the computation. We will then evaluate in \ref{sec:accuracy} with a number of test-case examples how the accuracy of a computed evolution can be ensured by selecting suitable values of the different accuracy parameters, including the integration stepsize, the mesh size and the remeshing frequency.

\subsection{Time stepping and stopping criteria}
\label{sec:timestep}

The overall or average curvature of the curve increases during an evolution as time $t$ converges to $t_{\max} < \infty$. This is why we determine the stepsize at step $i$ from the curvature as
\begin{equation}
h_i = \frac{S}{\max_{\varphi \in [-\pi,\pi]}(\kappa_i(\varphi))}.
\label{eq:timestepping}
\end{equation}
Here the accuracy parameter $S$ determines the fineness of the stepsize sequence $(h_i)$; it is chosen at the start of a computation and then remains constant. The computation stops when the time step $h_i$ drops below a pre-specified minimum $h_{\rm min}$.

During an evolution where the curve reaches a final point in $\bbR^2$ the radial distance function $r_i$ converges with $i$ to the constant function $0$; this also means that the area $A_i$ of $\Gamma_i$ goes to zero. Moreover, the average of the absolute value of the curvature $|\kappa_i|$ grows beyond bound. Working with too small $r_i$ or too large $|\kappa_i|$ is not practical numerically, which is why we limit their ranges; note that very small  $r_i$ or $A_i$ or very large curvature indicates that the total integration time $t_i = \sum_{l=1}^{i} h_l$ is close to $t_{\max}$. This is formalized in the following list of stopping criteria for the computation. \\[-3mm]
\begin{itemize}
\item[(i)] $h_i < h_{\rm min}$;
\item[(ii)] $\mathrm{Var}(r_{i+1}) < 1e^{-8}$;
\item[(iii)] $A_i < A_{min}$;
\item[(iv)] $\mathrm{max}_{\varphi \in [-\pi,\pi]}(|\kappa_{i}|/\mathrm{Avg}(r_{i}) < 1e^{-5}$. \\[-3mm]
\end{itemize}
Note that we impose global limits on the radial distance function $r_i$ and the curvature $\kappa_i$, while $h_{\rm min}$ and $A_{min}$ are accuracy parameters of the algorithm. The computation stops as soon as one of these stopping criteria is satisfied.

\subsection{Uniform meshes}
\label{sec:meshuni}

There are two types of uniform meshes that we consider as part of our implementation: meshes that are uniform in the phase $\varphi$, and meshes that are uniform in arclength $s$. We will use either as appropriate for the initial mesh $M_0$, as well as for remeshing during a computation. 

\begin{figure}[t!]
  \centering
  \includegraphics[scale=1]{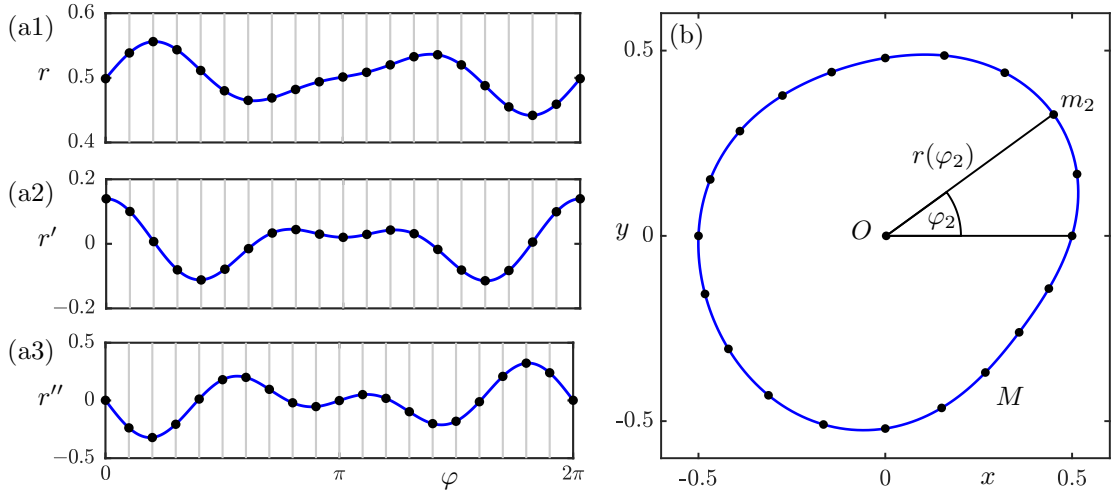}
  \caption{The piecewise polynomial central radial distance function $r$ (a1) and its derivatives $r'$ (a2) and  $r''$ (a3), generated by the phase-uniform mesh (black dots) of size $N=20$ on the standard example curve (b), as defined in \ref{eq:standardcurve}.}
  \label{fig:representation}
\end{figure}

A phase-uniform mesh $M$ for a star-like curve $\Gamma$ with polar parameterization $\gamma$ is given by
\begin{equation}
\label{eq:meshphi}
M^\varphi = \{\frac{2\pi}{N} j \ | \ 0 \leq j \leq N\} 
\quad {\rm with} \quad M = \gamma(M^\varphi).
\end{equation}
Phase-uniform meshes are computationally inexpensive to generate and their uniformity in $\varphi_j$ is an advantage for the construction of the piecewise polynomial periodic function $r$ of the associated discretization $\Pi_M(\Gamma) \in \mathscr{D}$. Figure \ref{fig:representation} shows the phase-uniform mesh of size $N = 20$ for our standard example curve, defined in \ref{eq:standardcurve} below, which generates the piecewise polynomial functions $r$, $r'$ and $r''$ of its discretization in $\mathscr{D}$. Note that the mesh points on $\Gamma$ are distributed almost uniformly in arclength as well. This is the case for any curve that is reasonably close to being circular, which is why, as we will see in \ref{sec:examples}, phase-uniform meshes are a suitable choice for such curves. 

An arclength-uniform mesh $M$ for $\Gamma$, on the other hand, is given by
\begin{equation}
\label{eq:mesharc}
M = \{\tilde{\gamma}(\frac{L}{N} j)  \ | \ 0 \leq j \leq N\} 
\quad {\rm with} \quad M^\varphi = \tilde{\gamma}^{-1}(M);
\end{equation}
here $\tilde{\gamma}$ is the associated arclength parameterization of $\gamma$ given by \ref{eq:polartoarc}. By construction, the mesh points $m_j$ are uniformly distributed in arclength along $\Gamma$, while the corresponding phases $\varphi_j \in M^\varphi$ are not uniformly distributed in $[-\pi, \pi]$. Working with arclength-uniform meshes is more computationally expensive but, as we will see in \ref{sec:examples}, may be required for the accurate computation of the evolutions of certain curves, in particular of ellipses. 

\begin{figure}[t!]
  \centering
  \includegraphics[scale=1]{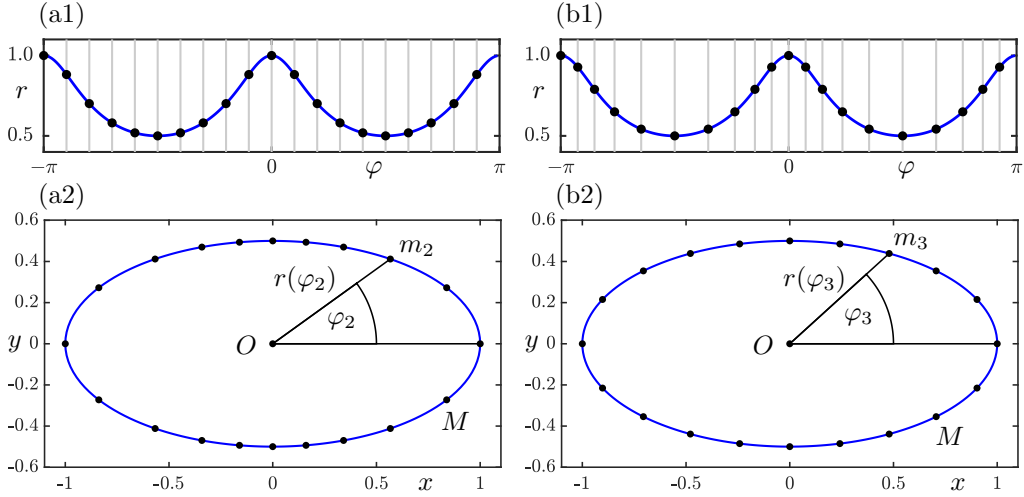}
  \caption{The piecewise polynomial central radial distance function $r$ generated by a uniform mesh of size $N=20$ on the standard ellipse \ref{eq:standardellipse} with $a = 1$ and $b = 0.5$, for a phase-uniform mesh in panels~(a1) and~(a2), and for an arclength-uniform mesh in panels~(b1) and~(b2).}
  \label{fig:discretizationellipse}
\end{figure}

Figure \ref{fig:discretizationellipse} illustrates the difference between the two types of uniform meshes with the example of a mesh of size $N = 20$ on an ellipse of axis ratio 2, the standard ellipse \ref{eq:standardellipse} with $a = 1$ and $b = 0.5$. The uniform mesh $M^\varphi$ generates a radial distance function $r$ in panel~(a1) with interpolated points that are at fixed distance $\frac{2\pi}{N}$. The associated meshpoints of the dual mesh $M$ on the ellipse in panel~(a2) are concentrated near the minima of the curvature and clearly not at equal arclengths from each other. For the arclength-uniform mesh the situation is reversed: the uniformity of the points of the mesh $M$ on the ellipse in panel~(b2) generates a non-uniformity of the associated interpolation points of $M^\varphi$ that define the radial distance function $r$ in panel~(b1). Observe how points in $M^\varphi$ now concentrate near phase angles corresponding to points of high curvature of the ellipse. Note that there is no distinguishable difference on the level of Figure \ref{fig:discretizationellipse} between the two radial distance functions and the two respective approximate ellipses; the accuracy of approximation for either uniform parameterization will be considered in more detail in \ref{sec:accuracy}.

\subsection{Mesh propagation}
\label{sec:meshprop}

Any computation of an evolution of any star-like curve $\Gamma$ starts with the choice of the initial mesh $M_0$ of size $N$ that then generates the piecewise polynomial central radial distance function $r_0$ that gives the parameterization of the initial curve $\Gamma_0 = \Pi_{M_0}(\Gamma) \in \mathscr{D}$. Throughout, we take $M_0$ to be a uniform mesh in either phase or arclength, depending on the situation; here $N$ is chosen sufficiently large to ensure that the errors between the respective periodic functions and their derivatives of $\Gamma$ and $\Gamma_0$ are sufficiently small; see \ref{sec:accuracy}. 

To compute the new radial distance function $r_{i}$ at step $i$ one needs to specify the new mesh $M_i$ on the curve $E_{h_i}(\gamma^{i-1})$. The most straightforward way to obtain $M_{i}$ is to push forward or propagate the mesh $M_{i-1}$ under the Euler step \ref{eq:Euler}. This means that $M_{i} = E_{h_i}(M_{i-1})$, which is computed as
\begin{equation}
\label{eq:meshEuler}
m_{i,j} = m_{i-1,j} + h_i (c + k_{i-1}^\alpha (\varphi_{i-1,j})) \n_{i-1,j}(\varphi_{i-1,j}));
\end{equation}
here $\varphi_{i-1,j} \in M^\varphi_{i-1} = \gamma^{-1}_{i-1}(M_{i-1})$ is the associated phase of the mesh point $m_{i-1,j} \in M_{i-1}$. 

\begin{figure}[t!]
  \centering
  \includegraphics[scale=1]{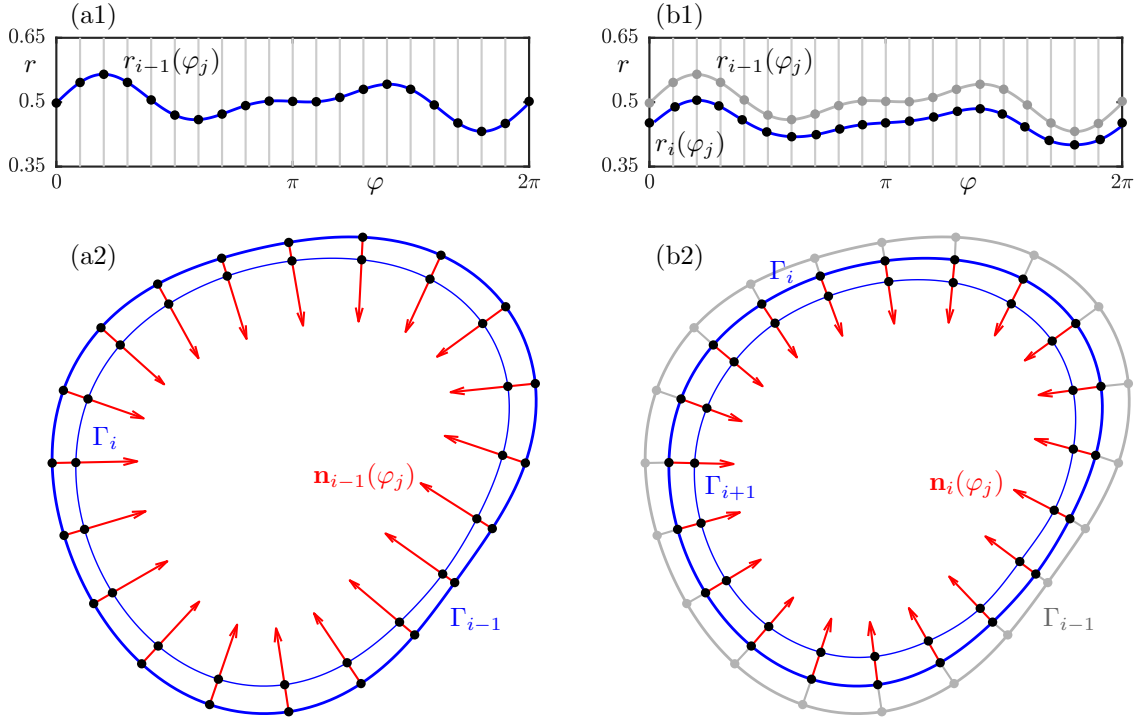}
  \caption{Illustration of step $i$ in panels (a) and of step $i+1$ in panels (b) of the mesh propagation algorithm for a mesh of size $N=20$. The radial distance function $r_{i-1}$ in panel~(a1) generates the parameterization of the curve $\Gamma_{i-1}$ in panel~(a2) with normals $\n_{i-1}$ (red arrows) at the mesh points (black dots). Euler steps at each mesh point give the new mesh that defines the next radial distance function $r_{i}$ shown in panel~(b1). It generates the parameterization of the next curve $\Gamma_{i}$  in panel~(b2) and its normals $\n_{i}$ (red arrows) at the mesh points (black dots). The figure is for the standard example curve from \ref{eq:standardcurve}; compare with Figure \ref{fig:representation}.} 
  \label{fig:algorithmsteps}
\end{figure}

This choice of propagating the mesh at each step gives the mesh propagation implementation of the general algorithm, which forms the core of our method. It has the advantage that all that is required computationally at each integration step is calculating the $N$ Euler steps at each mesh point, followed by the subsequent construction of the piecewise polynomial radial distance functions $r_i$. This is illustrated in Figure \ref{fig:algorithmsteps}, which shows steps $i$ and $i+1$ of the mesh propagation implementation.  The piecewise polynomial radial distance function $r_{i-1}$ in panel (a1) generates both the curve $\Gamma_{i-1}$ via its parameterization $\gamma_{i-1}(\varphi)$ as well as the normals $\n_{i-1}$ at the mesh $M_{i-1}$ in panel~(a2). The propagated mesh is the result of Euler steps at each mesh point, and it defines the piecewise polynomial radial distance function $r_{i}$ in Figure \ref{fig:algorithmsteps} (b1). The corresponding parameterization $\gamma_{i}(\varphi)$ and normals $\n_{i}$ at the mesh points in panel~(b2) then give the next curve with parameterization $\gamma_{i+1}(\varphi)$. Note that during each step of the computation there is a continuous switching between the two representations of the curve --- in the sense that we are using both the mesh of points on the respective curve, as well as the dual phase-mesh that defines the piecewise polynomial radial distance function needed to evaluate the normals and the curvature.

\subsection{Remeshing}
\label{sec:remesh}

During mesh evolution the mesh points tend to drift towards points of higher curvature of the curve. This effect is quite limited when the variation in curvature is small, but can be substantial when it is larger. Rather than countering this with an increase of the overall mesh size $N$, it is more efficient to perform remeshing at suitable steps of the evolution algorithm. This can be done in a straightforward way at step $i$ because the radial distance function $r_i$ and the parameterization $\gamma_i$ are functions that can be evaluated at any value of $\varphi \in [-\pi, \pi]$. Hence, remeshing simply consists of taking a uniform mesh $M^* \subset \Gamma_i$, either phase-uniform or arclength uniform, and computing the radial distance function defined by $M^*$. The result is the remeshed piecewise polynomial radial distance function $r_i^*$ with parameterization $\gamma^*_i$ of the reparameterized curve $\Gamma^*_i$. 

Remeshing requires the construction of the uniform mesh $M^*$ and the subsequent construction of the new radial distance function $r_i^*$, so adds to the computational effort. On the other hand, it allows one to work with meshes of smaller size $N$. There is clearly a balance in terms of computational effort and accuracy, when deciding on a remeshing frequency versus taking larger $N$. We also remark that it may be advantageous to increase or reduce the mesh size $N$ when remeshing. 

\begin{algorithm}[t!] 
\caption{Curve evolution} 
\label{tab:algorithm} 
\begin{algorithmic} 
	\REQUIRE a star-like closed curve $\Gamma$; parameters $c$ and $\alpha$ of \ref{eq:abflow}; number of meshpoints $N$; stepsize parameter $S$; accuracy parameters $h_{\max}$ and $A_{\rm min}$;\\[2mm]
    \ENSURE  construct the piecewise polynomial radial distance function $r_0$ for a uniform mesh $M_0$ of size $N$ with dual mesh $M^\varphi_0$ for a suitable reference point $O_0$; ensure that the associated curve $\Gamma_0$ with parameterization $\gamma_0$ is sufficiently close to $\Gamma$; set the step counter $i$ to 1;\\[2mm]
    \WHILE{none of criteria (i)--(iv) from \ref{sec:timestep} are satisfied}
		\STATE calculate the polynomial derivatives $r_{i-1}'$ 
and $r_{i-1}''$ symbolically;
		\STATE determine the stepsize $h_i$ from $\kappa_{i-1}$;
		\STATE evaluate the associated normals $\n_{i-1}$ and the 
curvature $\kappa_{i-1}$ in the points of $M^\varphi_{i-1}$;
		\STATE perform the Euler step \ref{eq:meshEuler} at the mesh points to obtain the propagated mesh \\ $M_{i} = E_{h_i}(M_{i-1})$ with dual $M^\varphi_i$;
		\STATE choose a suitable reference point $O_i$ for the 
new curve $\Gamma_{i} = E_{h_i}(\Gamma_{i-1})$;
		\STATE construct the piecewise polynomial radial distance function $r_i$ for $M_i$ and $O_i$;\\[1mm]
		\STATE \textbf{[when required:} remesh the curve with a uniform mesh $M^*$ by computing the associated piecewise polynomial radial distance function $r^*_i$ with respect to $O_i$; set $r_i = r^*_i$\textbf{]};\\[1mm]
                \STATE increase the step counter $i$ by 1;
    \ENDWHILE
\end{algorithmic}
\end{algorithm}

The overall evolution algorithm as implemented, with mesh propagation and remeshing when required, is presented in schematic form as Algorithm~\ref{tab:algorithm}. This formulation is for a general reference point $O_i$ at every step, but we stress once more that we are using the centroid $C_i$ as the reference point whenever it is available. 

\begin{figure}[t!]
  \centering
  \includegraphics[scale=1]{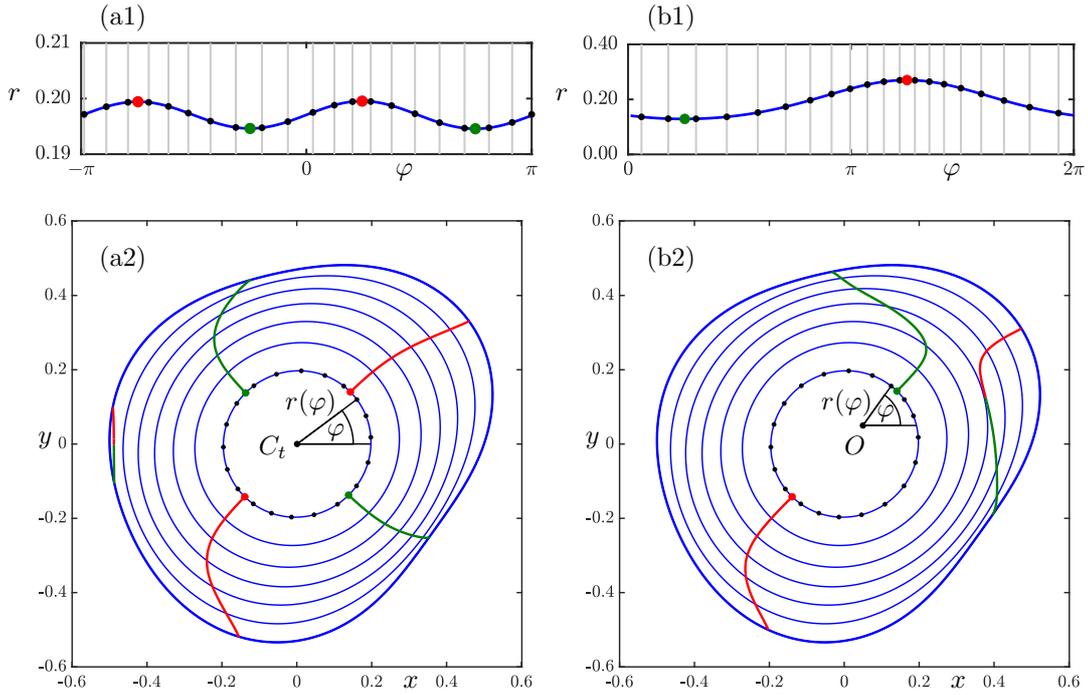}
  \caption{Evolution of the standard curve \ref{eq:standardcurve} under the curve-shortening flow, \ref{eq:abflow} with $c=0$ and $\alpha = 1$, as computed for stepsize parameter $S = 10^{-5}$ up to $t = 0.1062$  in 36,599 steps with a phase-uniform mesh of size $N =20$ at each step. Panels (a1) and (b1) show the radial distance functions of the final curve for the centroid $C_t$ and for $O = (0.05, 0.1)$ as the reference point, respectively. Panels (a2) and (b2) show seven curves at equally distributed time steps from $t=0$ to $t = 0.1062$, together with trajectories of the minima (green) and maxima (red) of the respective radial distance function as computed from the data at each step.}
  \label{fig:critpoints}
\end{figure}

A phase-uniform mesh $M^*$ is particularly efficient for remeshing because constructing it takes only $N$ evaluations of $r_i^*$; this means that the computational effort required at step $i$ is roughly doubled by the need to construct $r_i^*$. Phase-remeshing works well for curves that are close to being circular and evolve under a curvature flow that smoothes the curve. This is the case for the standard curve \ref{eq:standardcurve} from Figure \ref{fig:representation} under the curve-shortening flow. Its computed evolution up to $t_{\rm final} = 0.1062$ is shown in Figure \ref{fig:critpoints}; here a phase-uniform mesh of size $N =20$ was used with remeshing at every step, with stepsize parameter $S = 10^{-5}$ resulting in 36,599 integration steps. The calculation was performed with the centroid $C_i$ as the reference point; Figure \ref{fig:critpoints}(a1) and (a2) show the final central radial curve and the five  curves that are equally spaced in time between the initial curve and the final curve. Note that the final curve looks like a circle; however, panel (a1) shows that it has ondulations with four critical points, two minima and two maxima. Their images under the parameterization $\gamma_i$ are the stationary points of $\Gamma_i$, and they are shown at every step of the computation in Figure \ref{fig:critpoints}(a2). Panel (b1) shows the final radial curve for the off-centroid reference point $O = (0.05, 0.1)$. It has only two critical points, one minimum and one maximum. The images under $\gamma_i$ of the critical points of the non-central radial distance function at every step were computed by postprocessing and are shown in panel (b2) together with the same seven curves. Figure \ref{fig:critpoints} shows that this computation is accurate enough in either case to resolve the critical points of the radial distance function and, hence, the evolution and bifurcations of their images on the curve being evolved. Notice in panel (a2) how two stationary points disappear in a fold bifurcation early on, while the remaining four stationary points distribute more and more evenly around the curve as the evolution continues. In contrast, there are only four images of the critical points of the non-central radial distance function with reference point $O$, two of which disappear in a fold bifurcation while the remaining two distribute more and more evenly around the curve as the evolution continues. These observations are clearly in agreement with the statements in \ref{sec:known} regarding the numbers $n_{\mbox{\tiny $C$}}$ and $n_{\mbox{\tiny $O$}}$ and the evolution of the respective critical points.

Depending on the properties of the initial curve $\Gamma$ and the flow, a phase-uniform mesh may not be suitable because of a disadvantageous distribution of mesh points in arclength along the curve. In other words, a phase-uniform mesh would need to have an excessively large size in order to lead to a sufficiently accurate computation of the evolution. This is the case, in particular, for curvature flows with limits that are far from circular or develop corners; we will discuss such examples in \ref{sec:accuracy} and \ref{sec:andrewsalpha}. In such situation starting and remeshing with an arclength-uniform mesh is required. Constructing an arclength-uniform mesh is computationally more involved because it requires the calculation of the arclength integral in \ref{eq:polartoarc}; on the other hand, there may be considerable gain in working with a smaller mesh size.

In either case, the need to remesh depends on the quality of the mesh as determined by the arclength distribution along the curve. It is quite costly computationally to check at each step how even the arclength distribution of the mesh is; this is why we rather prescribe a frequency of remeshing suited to the evolution under consideration. For the evolution of the standard curve \ref{eq:standardcurve} in Figure \ref{fig:critpoints} we used a mesh of only $N=20$ and remeshed at each step, that is, with frequency 1. Depending on the starting curve, curvature flow, type of mesh and fineness of the time stepping, it is generally unnecessary to perform remeshing at every timestep. We will explore the effects of different remeshing strategies in more detail in \ref{sec:evolremesh}.

\section{Benchmark tests and accuracy settings}
\label{sec:accuracy}

We now consider how a computation is influenced by the different accuracy setting --- chiefly, the number of mesh point $N$, the type of mesh and the remeshing frequency. We first consider the interpolation error between a given star-like curve $\Gamma$ and its discretization in $\mathscr{D}$ and then compute evolutions for a number of benchmark test-case examples to determine suitable accuracy settings. In particular, we consider what we call the standard curve given by 
\begin{equation}
\label{eq:standardcurve}
r(\varphi) = 0.5 + 0.04\sin(2\varphi) + 0.03\sin(3\varphi),
\end{equation}
which was already seen in Figure \ref{fig:representation}. Moreover, we consider curves with additional symmetry, including the standard ellipse given by 
\begin{equation}
\label{eq:standardellipse}
r(\varphi) = \frac{ab}{\sqrt{(b\cos(\varphi))^2+(a\sin(\varphi))^2}}
\end{equation}
with $a > b > 0$; the ellipse of axis ration 2 shown in Figure \ref{fig:discretizationellipse} is given by \ref{eq:standardellipse} with $a = 1$ and $b = 0.5$.

\subsection{Projection onto initial piece-wise polynomial representation}
\label{sec:initialproj}

The first step of our algorithm is the discretization of the initial curve $\Gamma$ to obtain the initial representation $\Gamma_0 = \Pi_{M_0}(\Gamma) \in \mathscr{D}$. This introduces an interpolation error between $\Gamma$ and $\Gamma_0$ that can be measured as the distance between the respective radial distance functions and their derivatives, as defined by a norm on the space of periodic functions. Throughout, we consider the $L^2$ norm 
\begin{equation}
\label{eq:2normint}
\|v\|_2
= \left( \int_{-\pi}^{\pi} \| v(\varphi) \|^2 d \varphi 
\right)^{\frac{1}{2}} 
\end{equation}
of a periodic function $v$, as well as its associated discrete version
\begin{equation}
\label{eq:2normdiscr}
\|v\|^d_2
= \frac{1}{N} \left[ \sum_{j=1}^{N}  \| v(\varphi_j) \|^2 
\right]^{\frac{1}{2}} 
\end{equation}
that is evaluated only at the mesh points of $M_0$. 

\begin{figure}[t!]
  \centering
  \includegraphics[scale=1]{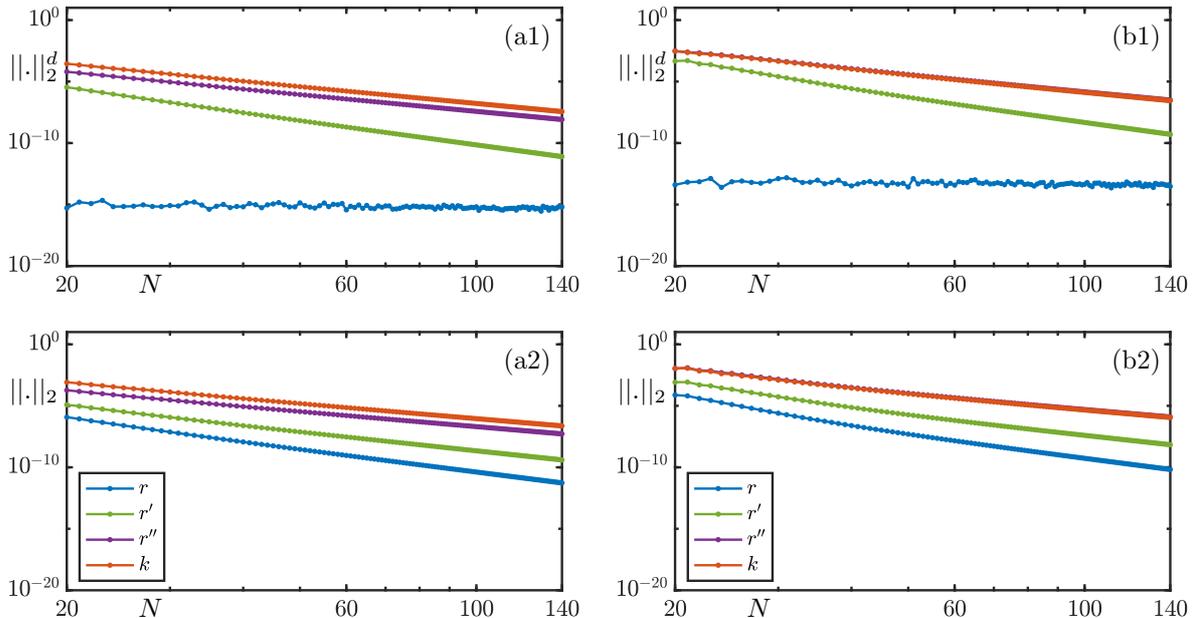}
  \caption{Discrete and integral error norms for the discretizations with a phase-uniform mesh of size $N$ with the centroid as reference point of the standard curve \ref{eq:standardcurve} in column (a) and of the standard ellipse \ref{eq:standardellipse} with $a = 1$ and $b = 0.5$ in column (b). 
}
  \label{fig:errornorms}
\end{figure}

To ensure the accuracy of the interpolation used in our implementation, we compute and monitor the integral and discretized error norms of the radial distance function, its first two derivatives and the signed curvature. These quantities are shown in Figure \ref{fig:errornorms} as a function of the mesh size $N$, for a phase-uniform mesh with the centroid as the reference point, of the standard curve \ref{eq:standardcurve} in column (a) and of the standard ellipse \ref{eq:standardellipse} with $a = 1$ and $b = 0.5$ in column (b). These doubly-logarithmic plots show that all errors converge to zero with the mesh size $N$, owing to the fact that the approximating curve $\Pi_{M_0}(\Gamma)$ converges to $\Gamma$ as expected from interpolation theory \cite{rivlin,trefethenbook}. Note that the discrete error in $r$ is at machine precision in both Figure \ref{fig:errornorms}(a1) and~(b1) since, by construction of the interpolating polynomials, the two curves agree at the mesh points. The interpolation error of $r$ in between mesh points is measured by the integral norm shown in panels~(a2) and~(b2). Also as expected, either error is larger for the derivatives $r'$ and $r''$, as well as for the derived quantity $\kappa$. (Note that $r''$ and $\kappa$ agree for the ellipse.) Overall, Figure \ref{fig:errornorms} shows that $N$ can be chosen sufficiently large to achieve a required accuracy of all quantities shown. Notice that to achieve a specified interpolation accuracy the mesh size $N$ needs to be chosen considerably larger for the ellipse compared to the standard curve. For example, the integral error in $\kappa$ for the standard curve with 60 mesh points is $7.3718e^{-6}$, and for the ellipse with 80 mesh points it is $2.0443e^{-5}$. This is due to the less optimal distribution of points in arclength along the ellipse; compare with Figure \ref{fig:representation} and Figure \ref{fig:discretizationellipse}.

\subsection{Mesh size during computed evolution}
\label{sec:evolmesh}

The interpolation error and, especially, the accuracy of $\n$ and $\kappa$ need to be controlled also during the computation of an evolution. In general, the evolution $\Gamma(t)$ of the given curve $\Gamma$ is not actually known. Therefore, the accuracy of a computation must be inferred by checking that, beyond some suitable mesh size $N$, the computed evolution is stable up to a desired accuracy under further increases of $N$. Since we are interested here in the associated evolution of the extrema of the radial distance function $r(\varphi;t)$ we now check, in particular, that their trajectories do not change with $N$, provided the mesh is fine enough. 

\begin{figure}[t!]
  \centering
  \includegraphics[scale=1]{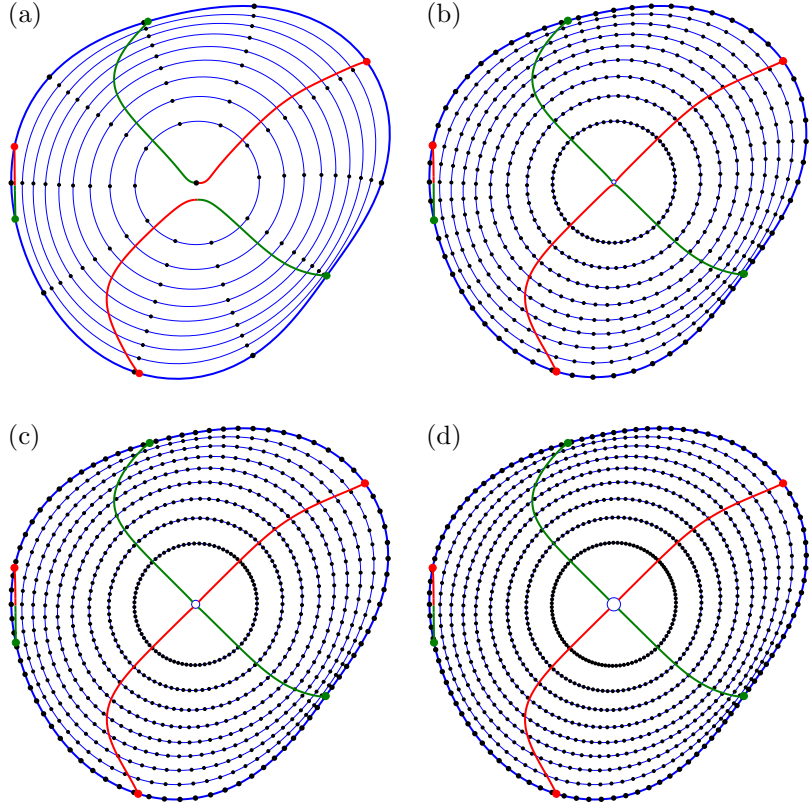}
  \caption{Evolution of the standard curve \ref{eq:standardcurve} under the curve-shortening flow, \ref{eq:abflow} with $c = 0$ and $\alpha = 1$, for a phase-uniform mesh with remeshing at every time step, for mesh size $N=10$ (a), $N=60$ (b), $N=80$ (c) and $N=100$ (d). Here $S = 10^{-5}$ and the reference point is the centroid; also shown are the trajectories of the minima (green) and maxima (red) of the central radial distance function.
}
  \label{fig:Nsize}
\end{figure}

Figure \ref{fig:Nsize} shows computations of the evolution of the standard curve \ref{eq:standardcurve} under the curve-shortening flow, given by \ref{eq:abflow} with $c = 0$ and $\alpha =1$, for phase-uniform meshes of four different sizes with remeshing at every time step. Here the time-step parameter is $S = 10^{-5}$, which lead to 56,519 steps in 89 s for $N=10$, to 57,272 steps in 391 s for $N=60$, to 58,239 steps in 515 s for $N=80$, and to 54,831 steps in 89s for $N=100$. These run times are for computations in Matlab on a standard desktop computer; the code has not been optimized and the run times are provided merely to give an indication of relative effort. In Figure \ref{fig:Nsize} and all similar figures only a few sample curves with meshes are shown, while the evolution of the trajectories of the extrema of the central radial distance function was derived from all computed time steps. On the level of the shown sample curves there is very little difference between the four panels of Figure \ref{fig:Nsize}. On the other hand, Figure \ref{fig:Nsize}(a) shows that a mesh of only $N=10$ points is definitely not sufficient for the resolution of the trajectories of the extrema. After all, they are the approximations of the stationary points and must, hence, reach the last computed curve (central point in (a), not shown). Close inspection of panel~(b) shows that for $N=60$ the lower two extrema meet just before the final curve is reached. In panels~(c) and~(d), on the other hand, the four remaining extrema all reach the central, final computed curve. Note that the final curve becomes larger with $N$, which is an indication that there are some issues with oversampling of a very small almost circular curve.

\subsection{Remeshing frequency}
\label{sec:evolremesh}

\begin{figure}[t!]
  \centering
  \includegraphics[scale=0.95]{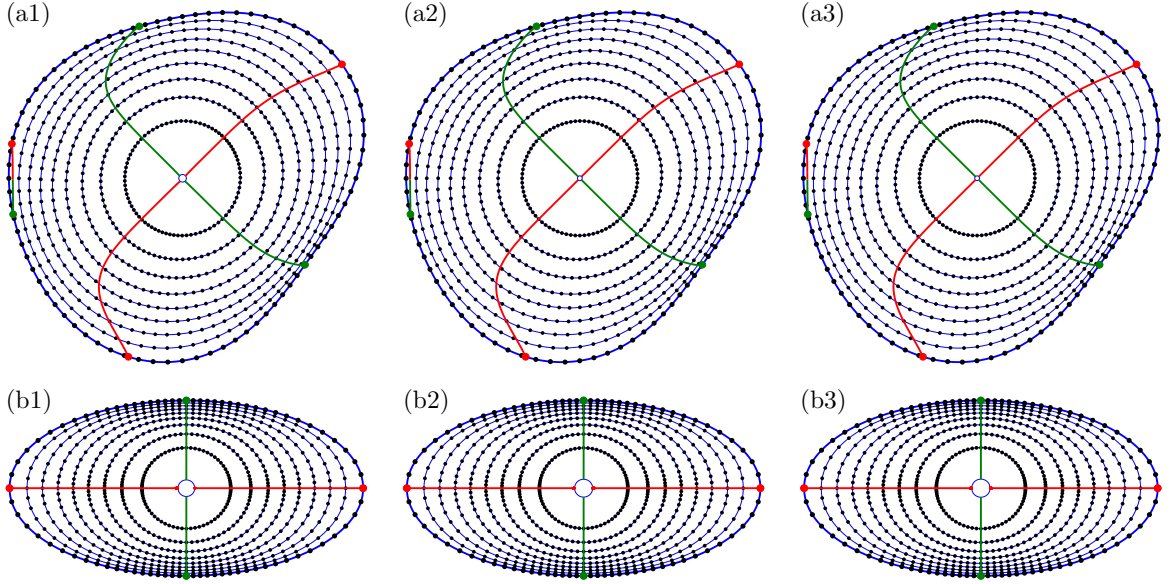}
  \caption{Evolution of the standard curve \ref{eq:standardcurve} in row (a) and of the ellipse \ref{eq:standardellipse} with $a = 1$ and $b = 0.5$ in row (b) under the
curve-shortening flow, \ref{eq:abflow} with $c = 0$ and $\alpha = 1$, computed with a phase-uniform mesh with $N=80$ points. Panels~(a1) and (b1) are computed without remeshing,  panels~(a2) and (b2) with remeshing at every 10th time step, and panels~(a3) and (b3) with remeshing at every time step. Here $S = 10^{-5}$ and the reference point is the centroid; also shown are the trajectories of the minima (green) and maxima (red) of the central radial distance function.
}
  \label{fig:remeshing_frequency}
\end{figure}

We settle on $N=80$ as a suitable mesh size of a phase-uniform mesh for the evolution of the standard curve as well as the ellipse with axis ratio 2 and now consider the effect of remeshing. Figure \ref{fig:remeshing_frequency} illustrates the respective computations of their evolution under the curve-shortening flow when there is no remeshing, that is, for propagation of the original mesh throughout, and for remeshing at every 10th and at every time step. The computations stop when the evolved curve is a circle within the numerical accuracy, which is the final circle in each of the panels of Figure \ref{fig:remeshing_frequency}. With a time-step parameter of $S = 10^{-5}$ this results for the standard curve in 55,541 steps in 382 s in (a1), 55,928 steps in 408 s in (a2) and 55,927 steps in 600 s in (a3); and for the ellipse in 110,991 steps in 433 s in (b1), 110,909 steps in 536 s in (b2) and 110,783 steps in 994 s in (b3). Note that, since both the standard curve and the ellipse with axis ratio 2 are not too far from being circular, the computations in panels~(a1) and (b1) without remeshing are already very accurate. In particular, the trajectories of the stationary points are resolved correctly; note that for the ellipse these are indeed horizontal and vertical lines. With remeshing at every 10th time step in panels~(a2) and (b2), one notices a slightly more uniform distribution of mesh points (in the phase variable $\varphi$), while remeshing at every time step, as in panels~(a2) and (b2), does not give any discernible improvement. 

\begin{figure}[t!]
  \centering
  \includegraphics[scale=0.95]{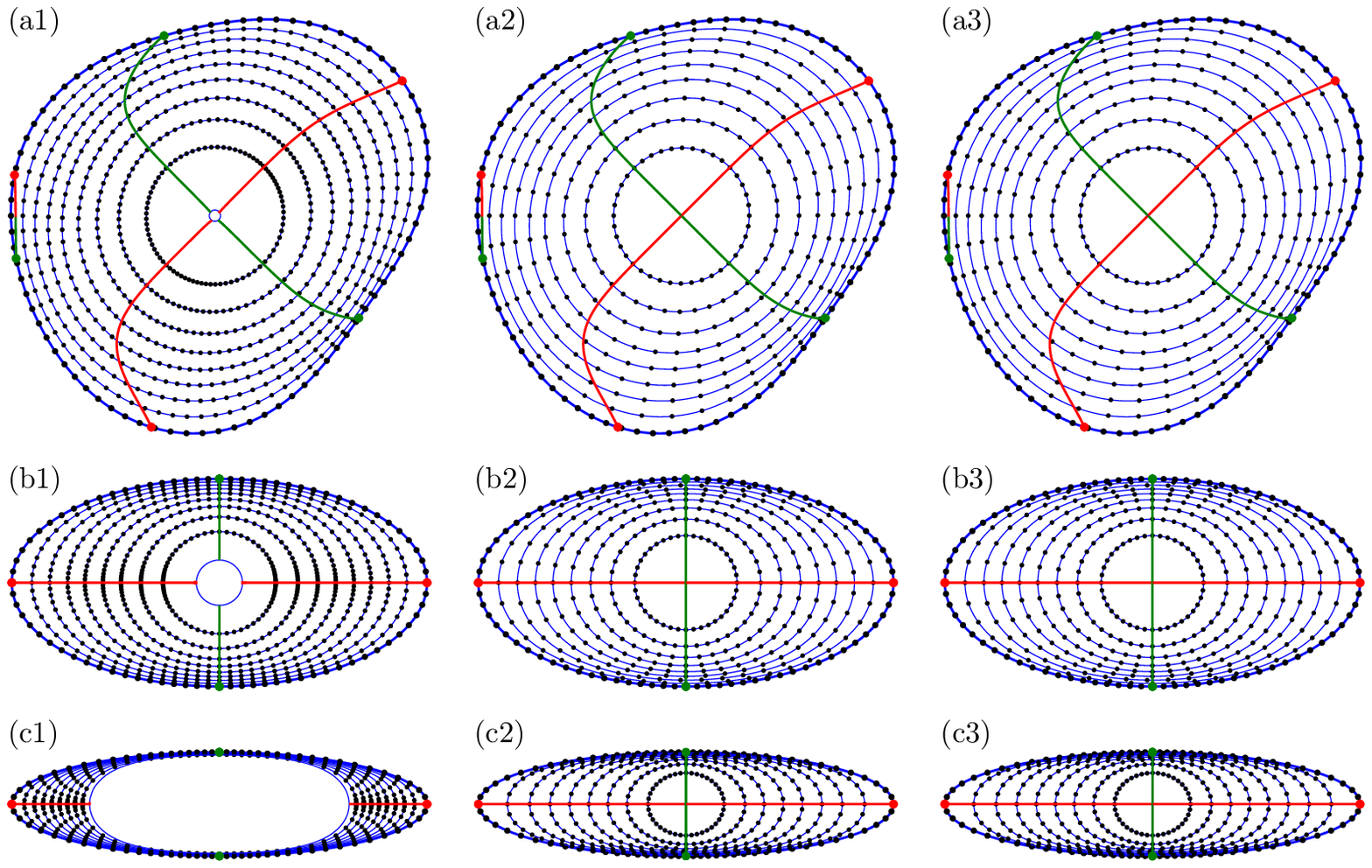}
\caption{Evolution of the standard curve \ref{eq:standardcurve} in row (a) and of the ellipse \ref{eq:standardellipse} with $a = 1$ and $b = 0.5$ in row (b) and of the ellipse \ref{eq:standardellipse} with $a = 2$ and $b = 0.5$ in row (c) under the curve-shortening flow, \ref{eq:abflow} with $c = 0$ and $\alpha = 1$, computed with an initial arclength-uniform mesh with $N=80$ points. Panels~(a1), (b1) and~(c1) are computed without remeshing,  panels~(a2), (b2) and~(c2) with remeshing at every 10th time step, and panels~(a3), (b3) and~(c3) with remeshing at every time step; during the remeshing the number $N$ of mesh points is reduced to keep the arclength distance $L/N$ between mesh points approximately constant, down to a minimum of $N_{min} = 40$ mesh points. Here $S = 10^{-5}$ and the reference point is the centroid; also shown are the trajectories of the minima (green) and maxima (red) of the central radial distance function. 
}
  \label{fig:remeshing_frequency_arclength_minN}
\end{figure}

We now consider the effect of remeshing when starting with a mesh of $N=80$ points that is uniform in arclength along the original curve. Here we again compute the evolutions of the standard curve and the standard ellipse of axis ratio 2, and also of an ellipse of axis ration 4, given by \ref{eq:standardellipse} with $a = 2$ and $b = 0.5$, which we refer to as the flat ellipse. In fact, the flat ellipse is quite far from being circular and computing it with a phase-uniform mesh would require a much larger number of mesh points. Figure \ref{fig:remeshing_frequency_arclength_minN} illustrates the respective computations of the three evolutions under the curve-shortening flow in rows~(a), (b) and~(c), in the absence of remeshing, with remeshing at every 10th and with remeshing at every time step. At each remeshing the number $N$ of mesh points is adjusted to keep the arclength distance $L/N$ between successive mesh points approximately constant; this leads to a reduction of $N$ down to a minimum of $N_{min} = 40$ mesh points. With $S = 10^{-5}$ this results for the standard curve in 55215 steps in 402 s in (a1), 56,431 steps in 1245 s in (a2) and 56,431 steps in 7,613 s in (a3); for the ellipse with axis ratio 2 in 89,398 steps in 567 s in (b1), 99,990 steps in 1,754 s in (b2) and 99,996 steps in 11,553 s in (b3); and for the flat ellipse with axis ratio 4 in 75,039 steps in 325 s in (c1), 199,469 steps in 4,214 s in (c2) and 199,452 steps in 19,795 s in (c3). 

Comparison of rows (a) and (b) of Figure \ref{fig:remeshing_frequency_arclength_minN} with those of Figure \ref{fig:remeshing_frequency} shows that the results are very similar: the trajectories are well resolved already without remeshing, remeshing at every 10th time step gives a better distribution of points along the computed curves and a closer approach to the final point, while remeshing at every time step does not give a noticeable further improvement. Remeshing of an arclength-uniform mesh is considerably more time consuming, compared to a phase-uniform mesh, owing to the need to evaluate the arlength integral \ref{eq:polartoarc}. On the other hand, it results in a very good mesh with a low number of mesh points. Row~(c) of Figure \ref{fig:remeshing_frequency_arclength_minN} shows that remeshing is crucial for curves that are further from circlular, that is, have segments of high curvature. Panel~(c1) for the flat ellipse shows that pure mesh propagation (no remeshing) leads to a strong accummulation of mesh points near the points of largest curvature; as a result, the computation stops very far from reaching the limit point. Remeshing clearly solves this problem, as panels~(c2) and~(c3) show; moreover, remeshing at every 10th time step again suffices. 

Our general algorithm with remeshing is flexible and efficient. The examples above show that regular remeshing allows one to work with relatively small meshes. On the other hand, remeshing requires additional computation time compared to pure mesh propagation, especially when one works with arclength-uniform meshes. We find for computed evolutions of curvature flows that become more circular that remeshing at every 10th time step provides a good balance between accuracy improvement and increased computational time. Moreover, for curves that are not too far from being circular, remeshing with phase-uniform meshes is accurate and especially fast. For curves with segments of large curvature, remeshing becomes a necessity and arclength-uniform meshes ensure accuracy with comparable numbers of mesh points. Note that for the chosen small value $S = 10^{-5}$ of the stepsize control parameter many thousands of curves $\Gamma_i$ are found as part of the computed evolutions, which allows us to compute the trajectories of critical points very accurately.

\subsection{Rotations and symmetry}
\label{sec:rotatedsymm}

Rotating a given curve $\Gamma$ in the plane results in the rotated evolution; in particular, the trajectories of the stationary points are identical up to the rotation. This property provides a good test for the accuracy of a computed evolution. In this context, the initial uniform mesh of a rotated curve generally has different mesh points on $\Gamma$; hence the piecewise polynomial radial distance function $r_0$ of the the projected initial curve $\Gamma_0$ changes under rotation. Indeed, for sufficiently stringent accuracy conditions, the computed evolution of the stationary points should be independent of the exact position of the initial mesh. This is indeed the case, as we show now.

\begin{figure}[t!]
  \centering
  \includegraphics[scale=1]{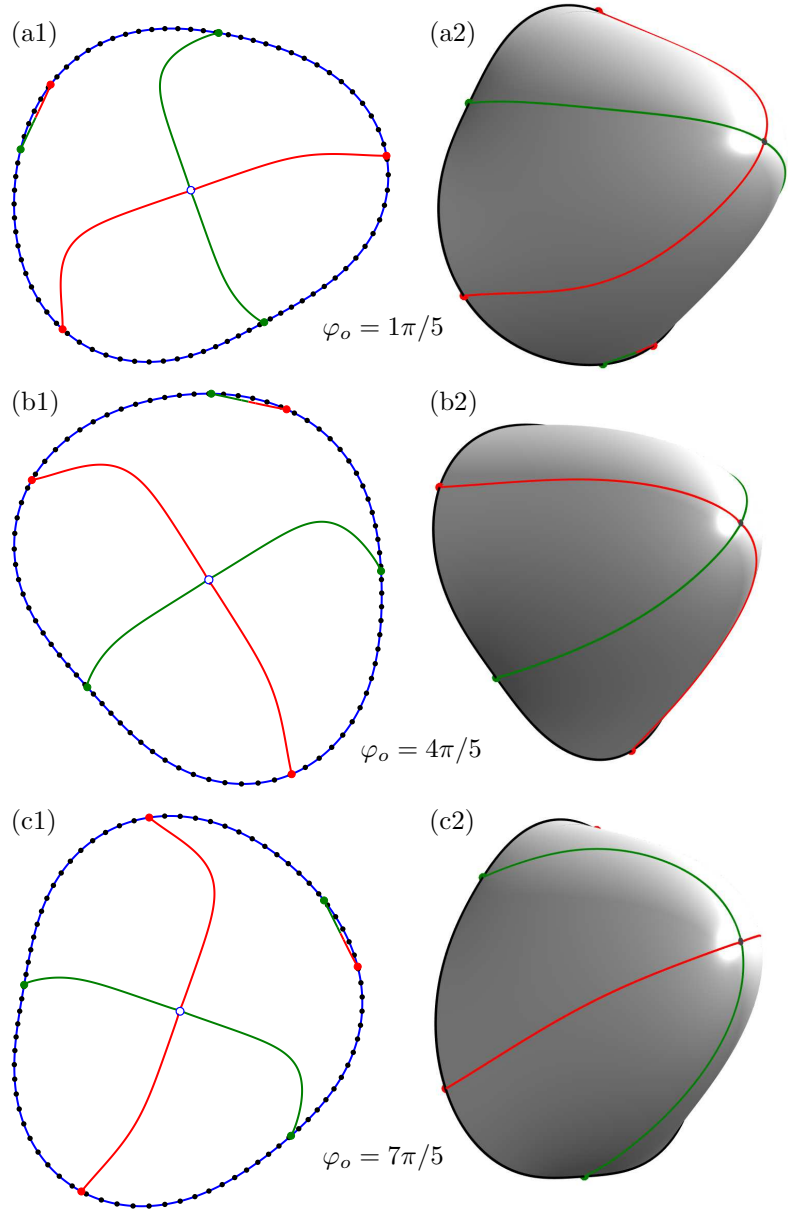}
  \caption{Evolution under the curve-shortening flow, \ref{eq:abflow} with $c = 0$ and $\alpha = 1$, of three rotations of the standard curve \ref{eq:standardcurve} over $\varphi_0$ as indicated, computed with a phase-uniform mesh of $N=80$ points with remeshing at every 10th time step. The left column shows the mesh on the standard curve and the trajectories of the stationary points, and the right column shows these trajectories on the evolution surface in $\bbR^2 \times \bbR^+$. Here $S = 10^{-5}$, the reference point is the centroid, and the stationary points are found as  the minima (green) and maxima (red) of the central radial distance function. 
}
  \label{fig:rotation}
\end{figure}

Figure \ref{fig:rotation} shows in rows (a)--(c) three evolutions under the curve-shortening flow of the standard curve \ref{eq:standardcurve} rotated over $\frac{1}{5}\pi$, $\frac{4}{5}\pi$ and $\frac{7}{5}\pi$, respectively. These evolutions were computed for a phase-uniform mesh with $N=80$ and remeshing at every 10th time step; the reference point is the centroid and $S = 10^{-5}$, which resulted in 55,926 steps in 359 s in (a), 55,901 steps in 379 s in (b) and 55,922 steps in 370 s in (c). The left column shows the trajectories of the stationary points, that is, the extrema of the central radial distance function, with only the initial curve $\Gamma_0$ with the phase-uniform mesh. The right column shows the trajectories of the stationary points on the evolution surface of the respective curve. This surface in Figure \ref{fig:rotation} has been rendered in $\bbR^3$, the product of $\bbR^2$ and the time axis $\bbR^+$, from the over 55,000 computed curves. The smoothness of this surface and of the trajectories of the stationary points on it illustrates the fineness of the time stepping and the overall accuracy of the computation in a new way. Note that all of these computations result in a very similar, but different sequence of computed curves $\Gamma_i$. Nevertheless, up to the respective rotation, the computed trajectories of the stationary points are indistinguishable from those in Figure \ref{fig:remeshing_frequency}(a2). 

\begin{figure}[t!]
  \centering
  \includegraphics[scale=1]{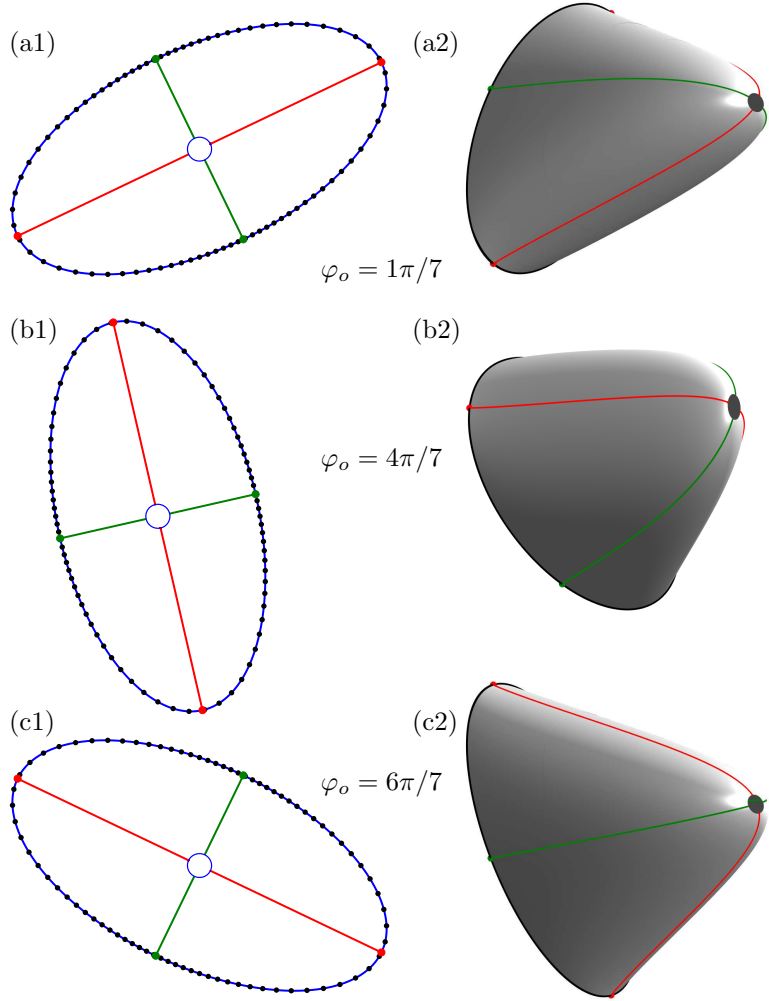}
\caption{Evolution under the curve-shortening flow, \ref{eq:abflow} with $c = 0$ and $\alpha = 1$, of three rotations of the standard ellipse \ref{eq:standardellipse} with $a=1$ and $b=0.5$ over $\varphi_0$ as indicated, computed with a phase-uniform mesh of $N=80$ points with remeshing at every 10th time step. The left column shows the mesh on the standard curve and the trajectories of the stationary points, and the right column shows these trajectories on the evolution surface in $\bbR^2 \times \bbR^+$. Here $S = 10^{-5}$, the reference point is the centroid, and the stationary points are found as the minima (green) and maxima (red) of the central radial distance function. 
}
  \label{fig:rotation_ellipse}
\end{figure}

Figure \ref{fig:rotation_ellipse} shows in the same way three evolutions under the curve-shortening flow of the standard ellipse \ref{eq:standardcurve} for $a = 1$ and $b = 0.5$ with axis ratio 2, rotated over $\frac{1}{7}\pi$ in row (a), $\frac{4}{7}\pi$ in row (b) and $\frac{6}{7}\pi$ in row (c). The computations are also for a phase-uniform mesh with $N=80$ with respect to the centroid, with remeshing at every 10th time step and time step control parameter $S = 10^{-5}$; this resulted in 99,235 steps in 591 s in (a), 99,288 steps in 542 s in (b) and 99,208 steps in 461 s in (c). The computations again generate slightly different steps and curves $\Gamma_i$ along the evolution and stop when the curve is the same effectively perfect circle; see the left column of Figure \ref{fig:rotation_ellipse}. The smoothness of the surfaces in the right column again illustrates the fineness of the time stepping. Note that the trajectories of the stationary points are perfectly aligned with the major and minor axes of the curves $\Gamma_i$, that is, are still straight lines in panels~(a1)--(c1), in spite of the fact that the extrema of the radial distance function occur in between mesh points. Also for the ellipse, the computed trajectories of the stationary points are practically indistinguishable from those in Figure \ref{fig:remeshing_frequency}(b2).  

Note, in particular, that in Figure \ref{fig:rotation_ellipse}, as well as in row~(b) of Figure \ref{fig:remeshing_frequency} and Figure \ref{fig:remeshing_frequency_arclength_minN}, the symmetry of the $D_2$-symmetry of the ellipse is preserved. Moreover, the trajectories of the stationary points are the fix-point subspaces of the reflections as required; see \ref{sec:symmGamma}. The preservation of the symmetry is another clear indication of the accuracy of our method.

\subsection{Homothetically contracting ellipses}
\label{sec:evolellipse}

In general, it is not known analytically how a given curve $\Gamma$ evolves under the flow of \ref{eq:abflow}. However, there are special (classes of) curves, called homothetically contracting solutions or homothetic curves, that do no change their shape under certain flows, meaning that $\Gamma(t)$ is a linear scaling of $\Gamma$ for all $0 < t < t_{\max}$. For some types of flows certain homothetically contracting solutions are known analytically. An immediate example are the circles, which do not change shape under \ref{eq:abflow} for any $0 < c$ and $0 < \alpha \leq 1$ \cite{Andrews2002a}. Note that every circle has a central parameterization $\gamma$ with radial distance function $r = const$, showing that the circles are in the class $\mathscr{D}$ of discretized curves we consider. Hence, our algorithm evolves circles as circles, without an approximation error. 

\begin{figure}[t!]
  \centering
  \includegraphics[scale=0.95]{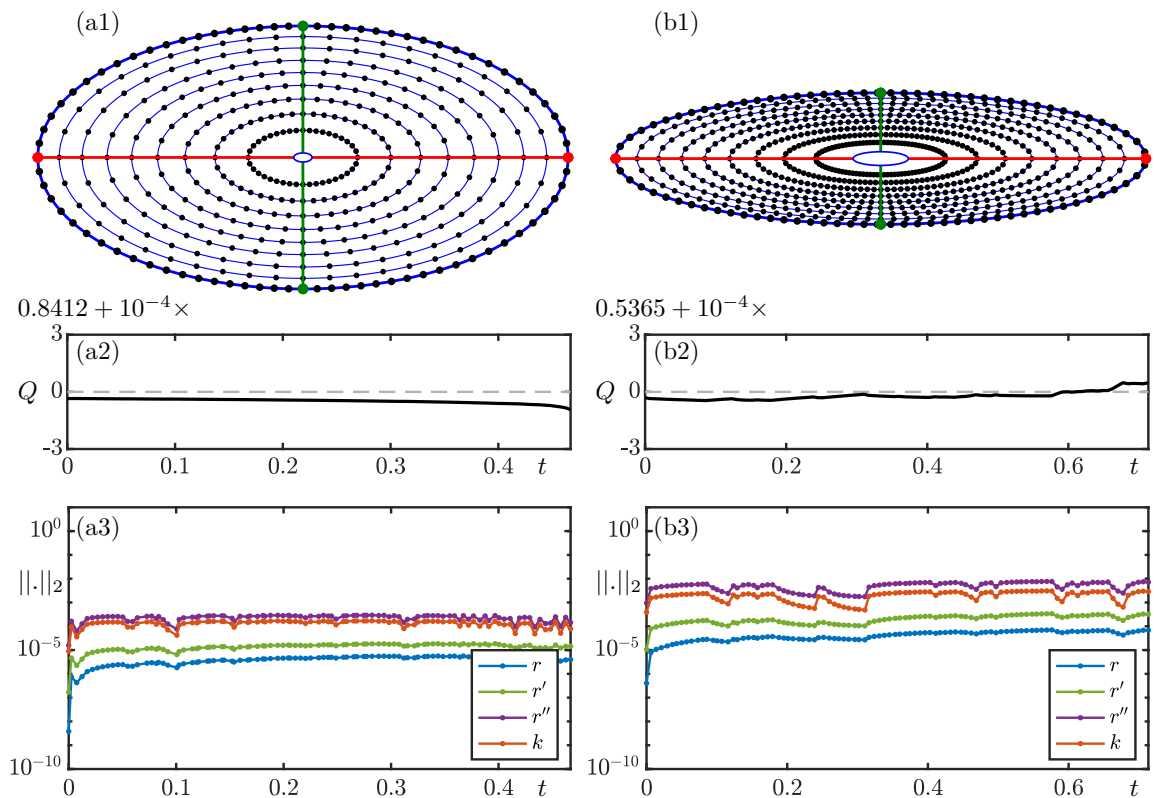}
  \caption{Evolution under the affine shortening flow, \ref{eq:abflow} with $c=0$ and $\alpha = \frac{1}{3}$, of the standard ellipse \ref{eq:standardellipse}, with $a=1$ and $b = 0.5$ in column (a) and with $a=2$ and $b = 0.5$ in column (b). Each computation is for an arclength-uniform mesh with $N=80$ points, $S = 10^{-5}$ and the centroid as the reference point. The trajectories of the stationary points in panels~(a1) and~(b1) show selected curves during the computation with the trajectories of the stationary points, found as the extrema of the central radial distance function. Panels~(a2) and~(b2) show the distance of the corresponding isoperimetric quotient $Q$ from its value for the respective ellipse as a function of time $t$, and panels~(a3) and~(b3) show the integral distance of the computed curves $\Gamma_{t_i}$ from the ellipse with same area and axis ratio, determined for the functions $r$, $r'$, $r''$ and $\kappa$.}
  \label{fig:erroraccumulation}
\end{figure}

According to Andrews \cite{Andrews1996,Andrews2002a}, a non-trivial example of analytically known homothetically contracting solutions are ellipses under \ref{eq:abflow} with $c=0$ and $\alpha = \frac{1}{3}$, which is known as the affine shortening flow. More specifically, starting from an initial ellipse, each curve $\Gamma(t)$ of this evolution is an ellipse with the same axis ratio. In particular, the isoperimetric quotient remains constant along the evolution. This special property of the affine shortening flow allows us to check the error along a computed evolution. Hence, we monitor during the computed evolution the isoperimetric quotient $Q$ (as computed for the polynomial representation from the area formula \ref{eq:computeA} and the arclength integral in \ref{eq:polartoarc}) as well as the distance of the (rescaled) curves $\Gamma_{t_i}$ to the initial ellipse.

Figure \ref{fig:erroraccumulation} shows in panels~(a1) and (b2) computed evolutions under the affine shortening flow of the ellipses from Figure \ref{fig:remeshing_frequency_arclength_minN}(b) and~(c) with axis ratio 2 and 4, respectively. Also shown in Figure \ref{fig:erroraccumulation} are the corresponding time series of the isoperimetric quotient $Q$ and of the integral distances between the curves $\Gamma_{t_i}$ and the respective ellipses. Throughout, we use an arclength-uniform mesh with $N=80$ points, remeshing every 10th time step, $S = 10^{-4}$ and the centroid as the reference point; this results in 70,886 steps in 1525 s in column~(a) and 227,646 steps in 5014 s in column~(b). Notice from panels~(a1) and~(b1) that the ellipses do not converge to a circle, but simply shrink while apparently maintaining their axis ratio. The latter is confirmed by the fact that $Q$ in panels~(a2) and~(b2) remains effectively constant during the evolution, staying within about $10^{-4}$ of its values 0.8412 and 0.5365 for the respective ellipses. Moreover, panels~(a3) and (b3) show that the integral error of the computation, computed as the integral distances of $r$, $r'$, $r''$ and $\kappa$ for $G_{t_i}$ and for the ellipse of the same area and with the same axis ratio, remains effectively constant as well. We conclude that, for the chosen mesh and accuracy settings, there is no significant error accumulation during the computation of either of these evolutions.

\section{Illustrations and numerical evidence for conjectures}
\label{sec:examples}

Our algorithm provides a new and efficient way to explore the solutions of the Andrews-Bloore flow \ref{eq:abflow} by performing numerical experiments. More specifically, it enables not only the computation of the evolution of an initial curve, but also allows us to follow accurately the trajectories of critical points of the radial distance function as well as other observables. We now demonstrate this capability. 

In \ref{sec:csf_trajectories} we consider the curve-shortening flow \ref{eq:curveshortening} and first illustrate known results regarding the evolutions of the numbers of critical points $n_{\mbox{\tiny $O$}}(t)$ and $n_{\mbox{\tiny $C$}}(t)$ with respect to a fixed reference point $O$ and the moving centroid $C(t)$, respectively; we then present numerical evidence for \ref{conj_U} on the properties of $n_{\mbox{\tiny $U$}}(t)$ with the fixed ultimate point $U$ as reference point and for \ref{conj_kappa} on the number of critical points $n_\kappa(t)$ of the curvature. Section \ref{sec:andrewsalpha} is devoted to illustrating results of Andrews regarding the shape evolution for $c=0$ and $0 < \alpha \leq 1$, and to finding the limit shapes of curves with $D_n$- and $C_n$-symmetry.

\subsection{Trajectories of critical points under the curve-shortening flow}
\label{sec:csf_trajectories}

We now illustrate known results discussed in  \ref{sec:known}, regarding the evolutions of the numbers $n_{\mbox{\tiny $O$}}(t)$ and $n_{\mbox{\tiny $C$}}(t)$ of critical points with respect to the fixed point $O$ and the moving centroid $C$, respectively, under the curve-shortening flow \ref{eq:curveshortening}, that is, under \ref{eq:abflow} with $c = 0$ and $\alpha = 1$. Figure \ref{fig:critpoints} already illustrated for the standard curve \ref{eq:standardcurve} that $n_{\mbox{\tiny $C$}}(t)$ decreases to 4, while the number $n_{\mbox{\tiny $O$}}(t)$ (with $O = (0.05,0.1)$) decreases to 2.

\begin{figure}[t!]
  \centering
  \includegraphics[scale=0.95]{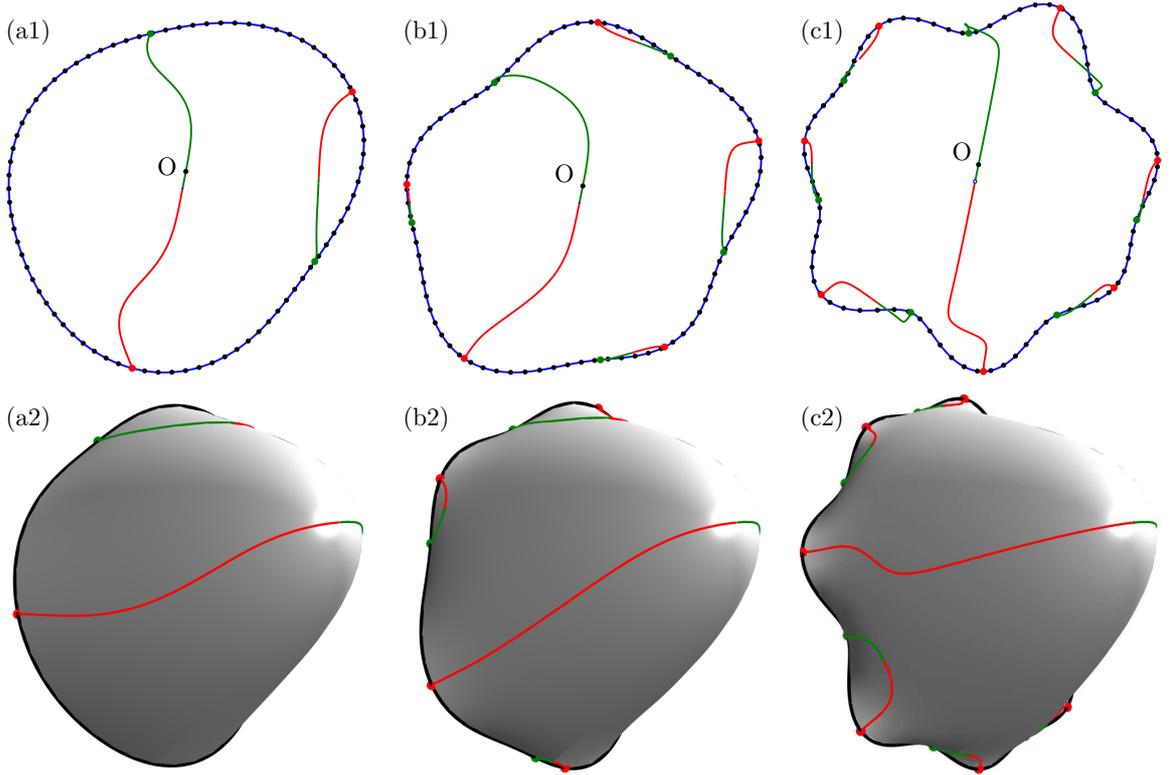}
  \caption{
Evolution under the curve-shortening flow, \ref{eq:abflow} with $c=0$ and $\alpha = 1$, of the trajectories of the critical points of the radial distance function $r_{\mbox{\tiny $O$}}$ for the fixed reference point $O = (0.01,0.05)$ (black dot) for the standard curve \ref{eq:standardcurve} in column (a), the curve given by \ref{eq:standardcurve5} in column (b), and the curve given by 
\ref{eq:standardcurve7} in column (c).  (Minima are shown in green and maxima in red.) Each computation is for an arclength-uniform mesh with $N=40$ points, $S = 10^{-5}$ and remeshing at every 10th time step. The top row shows the trajectories of the critical points of $r_{\mbox{\tiny $O$}}$ with the respective inital curve and mesh in $\bbR^2$, and the bottom row shows them on the evolution surface in $\bbR^2 \times \bbR^+$.}
  \label{fig:conjecture_01_05}
\end{figure}

\begin{figure}[t!]
  \centering
  \includegraphics[scale=0.95]{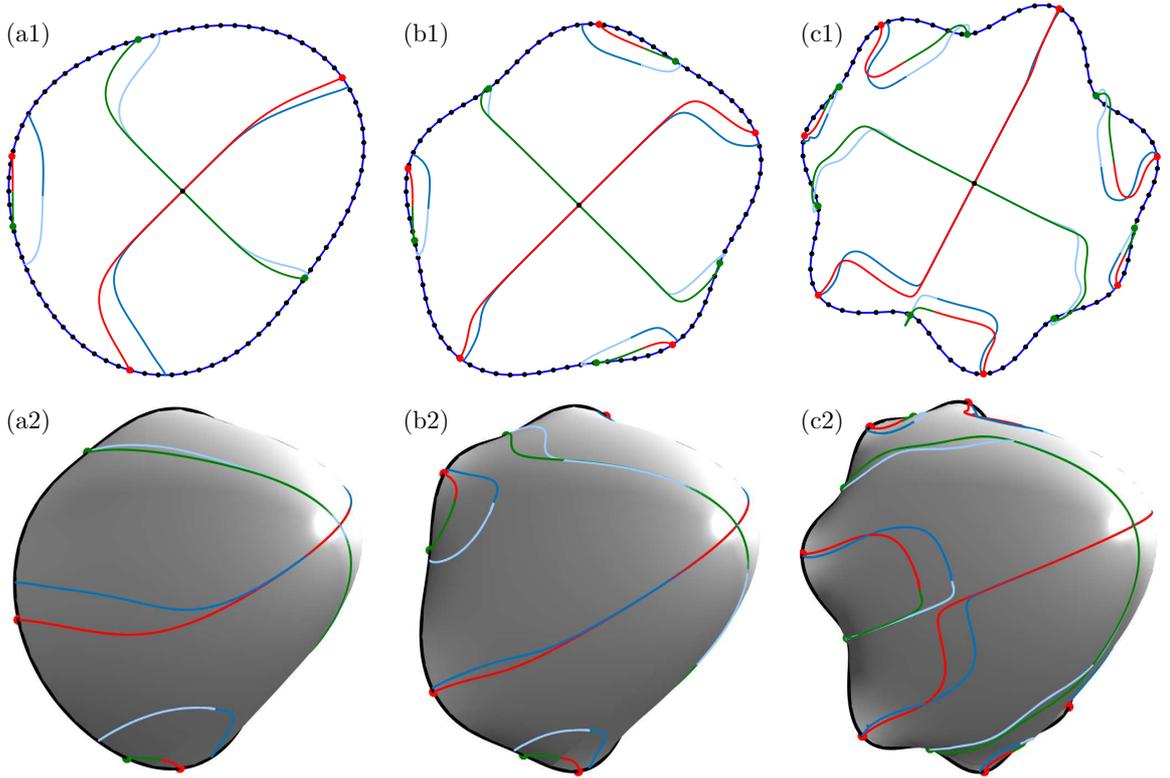}
\caption{
Evolution under the curve-shortening flow \ref{eq:curveshortening}, \ref{eq:abflow} with $c=0$ and $\alpha = 1$, of the trajectories of the critical points of the radial distance function $r_{\mbox{\tiny $C$}}$, with reference point at the (moving) centroid $C(t)$, for the standard curve \ref{eq:standardcurve} in column (a), the curve given by \ref{eq:standardcurve5} in column (b), and the curve given by \ref{eq:standardcurve7} in column (c).  (Minima are shown in green and maxima in red.) Also shown are the trajectories of the minima (light blue) and maxima (dark blue) of the signed curvature function $\kappa$.  Each computation is for an arclength-uniform mesh with $N=40$ points, $S = 10^{-5}$ and remeshing at every 10th time step. The top row shows the trajectories of the critical points of $r_{\mbox{\tiny $C$}}$ and of $\kappa$ with the respective initial curve and mesh in $\bbR^2$, and the bottom row shows them on the evolution surface in $\bbR^2 \times \bbR^+$.}
  \label{fig:conjecture_cxcy}
\end{figure}

We now consider the trajectories of the critical points of the radial distance function for three different reference points (fixed $O$, moving $C$ and fixed $U$) for three different curves. Namely, we compute the evolutions of
$n_{\mbox{\tiny $O$}}, n_{\mbox{\tiny $C$}}, n_{\mbox{\tiny $U$}}$ (and also that of $n_{\kappa}$) for the standard curve \ref{eq:standardcurve} as well as for the curves 
\begin{equation}
\label{eq:standardcurve5}
r(\varphi) = 0.5 + 0.04 \sin(2 \varphi) + 0.03 \sin(5 \varphi)
\end{equation}
and 
\begin{equation}
\label{eq:standardcurve7}
r(\varphi) = 0.5 + 0.03 \sin(5(\varphi+\pi/7)) + 0.04 \sin(7 \varphi) 
\end{equation}
with an arclength-uniform mesh of $N = 80$ points, $S = 10^{-5}$, and remeshing at every 10th time step. Any of the above-mentioned reference points could be chosen to perform the computation; we use the centroid $C$ throughout
and obtain the trajectories for the other three evolutions by postprocessing. 

Figure \ref{fig:conjecture_01_05} shows the trajectories of the critical points of the radial distance function $r_{\mbox{\tiny $O$}}(\varphi; t)$ for the fixed reference point $O = (0.01,0.05)$. Their number $n_{\mbox{\tiny $O$}}(t)$ decreases monotonically to $n_{\mbox{\tiny $O$}}(t_{max}) = 2$ at generic fold bifurcations. We expect to see more degenerate bifurcation if $\Gamma$ has some nontrivial symmetries; see the examples in \ref{sec:evolution}. 
These findings are in agreement with the statement proven in \cite{Domokos2015} that for $0 < \alpha$ the number $n_{\mbox{\tiny $O$}}(t)$ (for a reference point) decreases monotonically. Moreover, Figure \ref{fig:conjecture_01_05} illustrates the observation made in \ref{sec:known} that $n_{\mbox{\tiny $O$}}(t_{\max})=2$ and the two trajectories of critical points meet at 180 degrees (when continued to the ultimate point $U$).

Figure \ref{fig:conjecture_cxcy} shows the trajectories of the critical points of the central radial distance function $r_{\mbox{\tiny $C$}}(\varphi; t)$ with the reference point at the (moving) centroid $C$, along with the trajectories of the critical points (minima: light blue, maxima: dark blue) of the signed curvature $\kappa(\varphi; t)$. As can be seen, $n_{\mbox{\tiny $C$}}(t)$ monotonically decreases at successive fold bifurcations to $n_{\mbox{\tiny $C$}}(t_{max}) = 4$ and, similarly, $n_{\kappa}$ decreases until $n_{\kappa}(t_{max}) = 4$. In both cases the four trajectories meet at 90 degrees with one another. Notice also that $n_{\mbox{\tiny $C$}}(t_{max}) = 4$ and $n_{\kappa}$ are quite closely aligned with one another as the limit $t_{max}$ is approached. These findings support the statement regarding $n_{\mbox{\tiny $C$}}(t)$ in \ref{sec:known} as well as \ref{conj_kappa}. Moreover, they indicate that, as $t \to t_{\max}$, the curve $\Gamma(t)$ approaches a curve with $D_2$-symmetry, further confirming Bloore's result about the decay of Fourier terms \cite{Bloore1977}.

For any of the three initial curves shown in Figure \ref{fig:conjecture_cxcy}, the trajectories of the critical points of the radial distance function $r_{\mbox{\tiny $U$}} (t)$, with the ultimate point $U$ as the fixed reference point, is numerically indistinguishable from the evolution of the critical points of $r_{\mbox{\tiny $C$}}(\varphi; t)$. This means, in particular, that $n_{\mbox{\tiny $U$}}(t)$ decreases monotonically to $n_{\mbox{\tiny $U$}}(t_{max}) = 4$ at generic fold bifurcations and that the trajectories meet at 90 degrees, which is evidence for \ref{conj_U}. The effectively identical nature of the critical points of $r_{\mbox{\tiny $U$}}(\varphi; t)$ and of $r_{\mbox{\tiny $C$}}(\varphi; t)$ can be attributed to the fact that the centroid $C(t)$ practically does not move during the respective evolutions. We remark that it is an interesting task beyond the scope of this paper to find an initial curve with a significant difference between $n_{\mbox{\tiny $U$}}(t)$ and $n_{\mbox{\tiny $C$}}(t)$.

\subsection{Shape evolution for $c = 0$ and $0 < \alpha \leq 1$}
\label{sec:andrewsalpha}

We now illustrate some intriguing properties of the Andrews flow, that is, of \ref{eq:abflow} for $c=0$ (that is, without a constant term). As was mentioned already in \ref{sec:known}, following work by Gage~\cite{Gage1984,Gage1986} and Grayson~\cite{Grayson1987} for the case $\alpha  = 1$, Andrews examined how the limit of a curve under evolution depends on $0 < \alpha$ \cite{Andrews1996,Andrews1998,Andrews2002a}. A main result of his is the dichotomy with respect to the special case $\alpha = \frac{1}{3}$ of the affine shortening flow: for $\frac{1}{3} < \alpha$ any smooth convex curve converges to a round point and  the isoperimetric quotient converges to 1, while for $0 < \alpha \leq \frac{1}{3}$ the circle is not a limit any more. More specifically, the statement is that for any $j \geq 0$ there exists an open and dense set in the space of convex $C^j$-smooth curves such that the evolution of any curve from this set has isoperimetric quotient approaching zero as $t \to t_{\max}$ \cite[Theorem~7.1]{Andrews2002a}; this is an improvement of an earlier result in \cite{Andrews2002b} that required some symmetry of the curves.
For the dividing case $\alpha = \frac{1}{3}$ of the affine shortening flow any ellipse is a homothetically contracting solutions, which is the property we used for our accuracy test in \ref{sec:evolellipse}.

\begin{figure}[t!]
  \centering
  \includegraphics[scale=0.95]{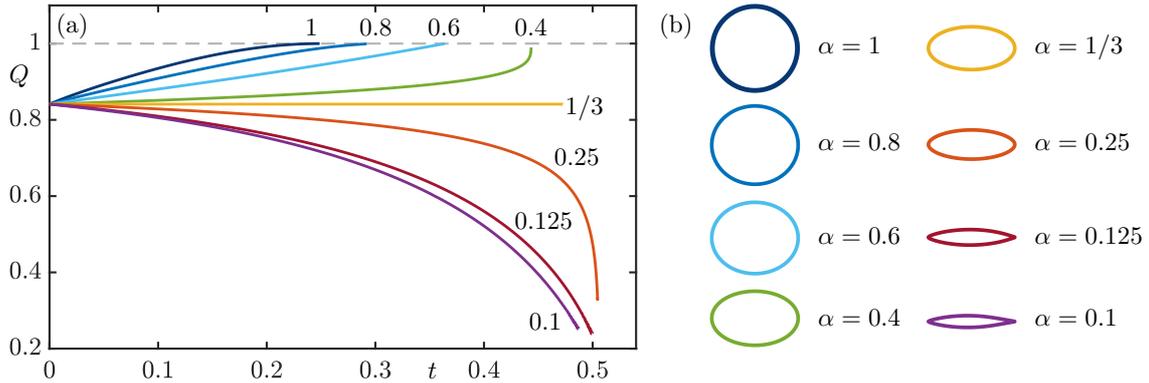}
\caption{Evolution of the standard ellipse \ref{eq:standardellipse} with $a=1$ and $b=0.5$ under~\ref{eq:abflow} with $c = 0$ for $\alpha = 1, 0.8,  0.6, 0.4, \frac{1}{3},  0.25,  0.125, 0.1$, computed with an arclength-uniform mesh with $N=80$, remeshing at every 10th time step and $S=10^{-5}$. Panel~(a) shows the respective time series of the isoperimetric quotient $Q$ and panel~(b) presents the rescaled final curves.}
  \label{fig:isoalpha}
\end{figure}

In other words, for any typical curve the isoperimetric quotient $Q$ will converge to 1 when $\frac{1}{3} < \alpha < 1$ and to 0 when $0 < \alpha < \frac{1}{3}$; moreover, for the special case of ellipses, $Q$ remains constant when $\alpha = \frac{1}{3}$. To illustrate this dichotomy, Figure \ref{fig:isoalpha} presents the evolution of the standard ellipse \ref{eq:standardellipse} with $a = 1$ and $b = 0.5$  under the flow of~\ref{eq:abflow} for $c=0$ for eight values of $\alpha$; here panel~(a) shows the respective time evolutions of the isoperimetric quotient $Q$ and panel (b) shows the corresponding final curves. The computations are for an arclength-uniform mesh with $N=80$ points and remeshing at every 10th time step. In accordance with \ref{sec:evolellipse}, Figure \ref{fig:isoalpha}(a) shows a constant $Q \approx 0.8412$ for $\alpha = \frac{1}{3}$ and an unchanged shape of the ellipse in panel~(b). Moreover, $Q$ indeed converges to 1 for $\frac{1}{3} < \alpha$, but note that the convergence of $Q$ to 1 and, hence, that of the curve to the circle, becomes considerably slower as $\alpha$ is decreased from 1. Similarly, Figure \ref{fig:isoalpha}(a) shows that $Q$ indeed decreases towards $0$ for $0 < \alpha < \frac{1}{3}$; again, this decrease is slower the closer $\alpha$ is to $\alpha = \frac{1}{3}$. The computed final curves in panel~(b) show that the corresponding curves indeed no longer approach a circular shape but, rather, flatten out. In the process, the curvature at the two points of maximal curvature on the ellipse appears to increase beyond any bound, while the curvature along two arcs connecting them decreases towards zero. We conclude that the curve approaches a double-covered straight line segment and the radial distance function $r_{\mbox{\tiny $C$}}$, when rescaled to fixed arclength $L = 4$ of the curve, converges to the discontinuous function that takes the value $1$ at $0, \pi$ and is identically zero elsewhere. This explains why the limit cannot be reached with our algorithm, because it assumes smoothness of $r_{\mbox{\tiny $C$}}$. Nevertheless, Figure \ref{fig:isoalpha} is clear numerical evidence that for $0 < \alpha < \frac{1}{3}$ the isoperimetric quotient $Q$ approaches zero as  $t \to t_{\max}$. 

 \begin{figure}[t!]
   \centering
\includegraphics[scale=0.95]{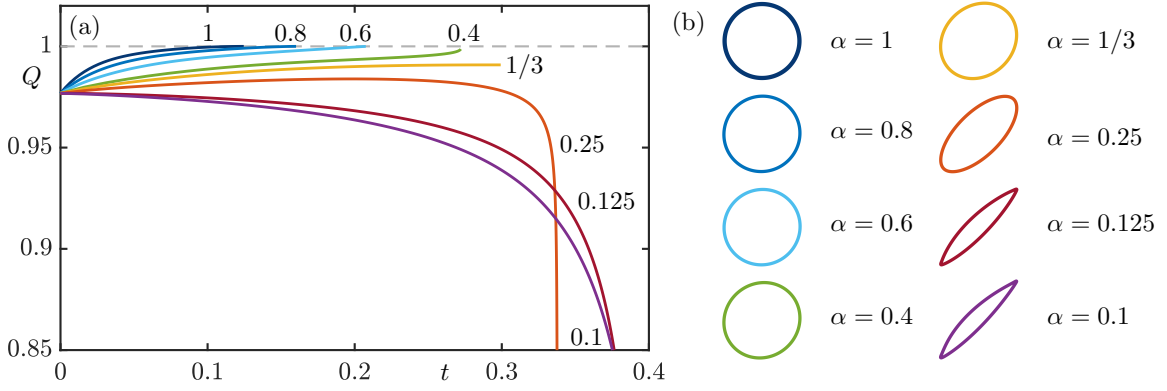}
   \caption{Evolution of the standard curve \ref{eq:standardcurve} under \ref{eq:abflow} with $c=0$ for $\alpha = 1, 0.8,  0.6, 0.4, \frac{1}{3},  0.25,  0.125, 0.1$, computed with an arclength-uniform mesh with $N=80$, remeshing at every 10th time step and $S=10^{-5}$. Panel~(a) shows the respective time series of the isoperimetric quotient $Q$ and panel~(b) presents the rescaled final curves.}
   \label{fig:standard_isoperimetric}
 \end{figure}

We now show that the same phenomenon occurs for a typical curve without any symmetry. Figure \ref{fig:standard_isoperimetric} shows the evolution of the standard curve \ref{eq:standardcurve} under \ref{eq:abflow} with $c=0$ for the same values of $0 < \alpha \leq 1$; again, panel~(a) presents the time series of $Q$, while panel (b) shows the corresponding final curves. The computations are for an arclength-uniform mesh with $N=80$ points and remeshing at every 10th time step with $S=10^{-5}$. Also for the standard curve there is a clear dichotomy: for $\frac{1}{3} < \alpha \leq 1$ the curve becomes a circle in the limit with $Q$ converging to $1$, while for $0 < \alpha \leq \frac{1}{3}$ the isoperimetric quotient decreases towards $0$. Note that in Figure \ref{fig:standard_isoperimetric}(a) we only show $Q$ over the range $[0.85, 1]$; for $\alpha = 0.25,  0.125, 0.1$ the computation reaches the values $Q = 0.788, 0.386, 0.337$, respectively. While these evolutions stop well before reaching $Q = 0$, Figure \ref{fig:standard_isoperimetric} nevertheless constitutes numerical evidence that also a general curve attains this limit by converging to a double-covered straight interval (after rescaling). In particular, the radial distance function $r_{\mbox{\tiny $C$}}$, when rescaled to fixed arclength $L = 4$ of the curve, converges to the discontinuous function that takes the value $1$ at $\varphi^*, \varphi^*+\pi$ and is identically zero elsewhere. Notice further from Figure \ref{fig:standard_isoperimetric} that the evolution for $\alpha = \frac{1}{3}$ converges to an ellipse, with $Q$ reaching a corresponding constant value; this illustrates the statement from \cite{Andrews1996} that any convex curve converges to an ellipse for $\alpha = \frac{1}{3}$. It appears that the (rescaled) interval on which the curve converges for $0 < \alpha \leq \frac{1}{3}$ is aligned with the long axis of this ellipse, which is in agreement with what we found for the standard ellipse in Figure \ref{fig:isoalpha}. 

\subsubsection{Curves with higher symmetry}
\label{sec:highersymm}

\begin{figure}[t!]
  \centering
  \includegraphics[scale=1]{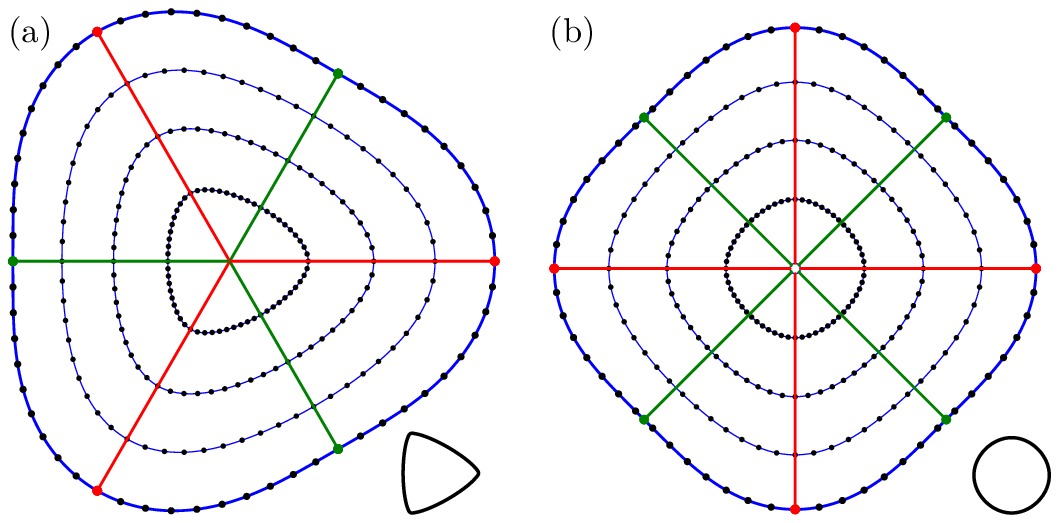}
  \caption{Evolutions under \ref{eq:abflow} with $c=0$ and $\alpha = 0.1$ of the $D_n$-symmetric curves \ref{eq:standardsymmcos} with $n = 3$ and $\delta = 0.05$ (a), and with $n = 4$ and $\delta = 0.03$ (b), computed with an arclength-uniform mesh with $N=60$ for $n =3$ and $N = 64$ for $n = 4$, remeshing at every 10th time step and $S=10^{-4}$. Shown are a selection of curves with their respective meshes and the trajectories of the stationary points; the computed final shape of the respective evolution is shown at the bottom right of each panel.
}
  \label{fig:polygonal_kappa01}
\end{figure}

\begin{figure}[t!]
  \centering
  \includegraphics[scale=1]{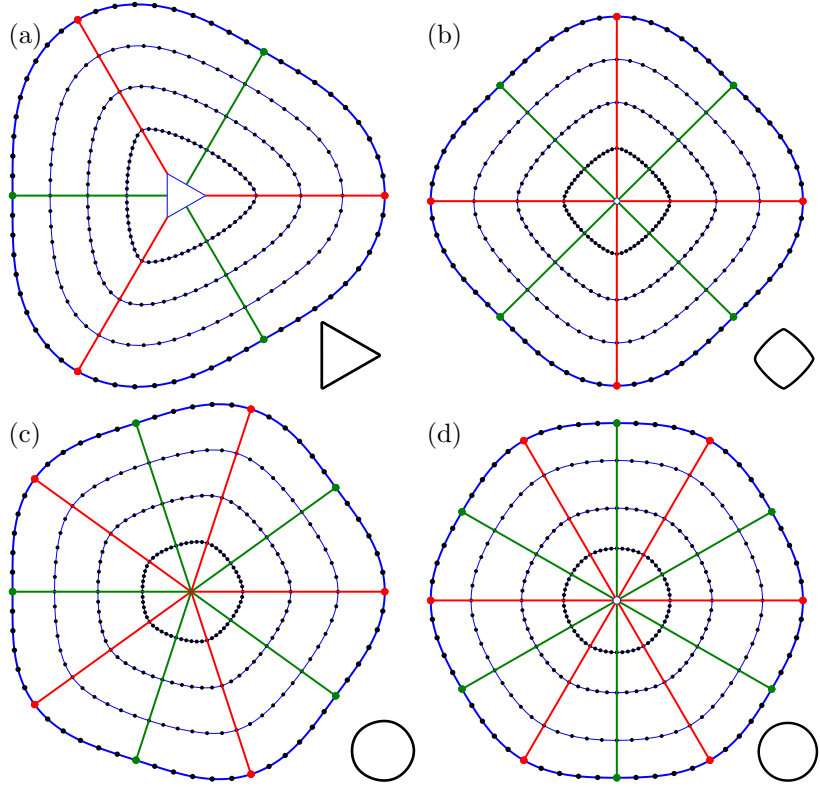}
  \caption{Evolutions under \ref{eq:abflow} with $c=0$ and $\alpha = 0.05$ of the $D_n$-symmetric curves \ref{eq:standardsymmcos} with $n = 3$ and $\delta = 0.05$ (a), with $n = 4$ and $\delta = 0.03$ (b), with $n = 5$ and $\delta = 0.02$ (c), and with $n = 6$ and $\delta = 0.01$ (d), computed with an arclength-uniform mesh with $N = 64$ for $n = 4$ and  $N=60$ otherwise, remeshing at every 10th time step and $S=10^{-4}$. Shown are a selection of curves with their respective meshes and the trajectories of the stationary points; the computed final shape of the respective evolution is shown at the bottom right of each panel.
}
  \label{fig:polygonal_kappa005}
\end{figure}

\begin{figure}[t!]
  \centering
  \includegraphics[scale=1]{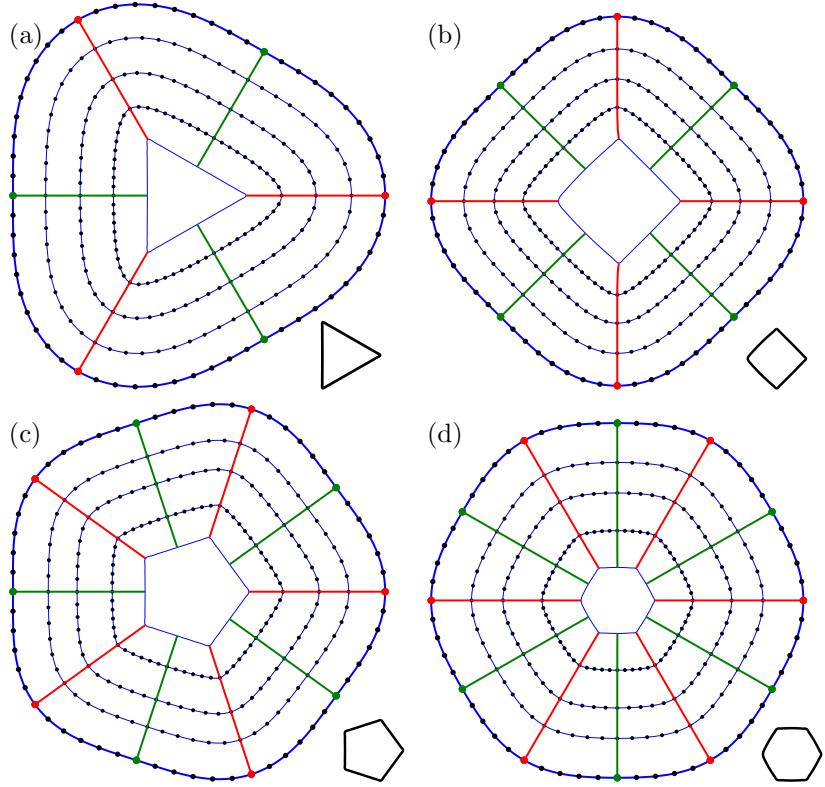}
  \caption{Evolutions under \ref{eq:abflow} with $c=0$ and $\alpha = 0.01$ of the $D_n$-symmetric curves \ref{eq:standardsymmcos} with $n = 3$ and $\delta = 0.01$ (a), with $n = 4$ and $\delta = 0.03$ (b), with $n = 5$ and $\delta = 0.02$ (c), and with $n = 6$ and $\delta = 0.01$ (d), computed with an arclength-uniform mesh with $N = 64$ for $n = 4$ and $N=60$ otherwise, remeshing at every 10th time step and $S=10^{-4}$ and the centroid as the reference point. Shown are a selection of curves with their respective meshes and the trajectories of the stationary points; the computed final shape of the respective evolution is shown at the bottom right of each panel.
}
  \label{fig:polygonal_kappa001}
\end{figure}

Andrews \cite{Andrews2002} provided further insight into \ref{eq:abflow} with $c=0$. First of all, for any $\alpha$, he showed that circles are homothetically contracting solutions and they are the \emph{only} homothetic solutions for $\frac{1}{8} \leq \alpha < \frac{1}{3}$ and $\frac{1}{3} < \alpha \leq 1$. On the other hand, for $0 < \alpha < \frac{1}{8}$ there are, additionally, unique homothetic embedded curves (up to translation, rotation and scaling) with $D_n$-symmetry. For any $0 < \alpha < \frac{1}{8}$ and for every $n$ satisfying $n < \sqrt{\frac{\alpha+1}{\alpha}}$ there exists such a unique embedded curve. Expressed in terms of $\alpha$, this condition becomes 
\begin{equation}
\label{eq:dncondition} 
\alpha < \frac{1}{n^2 - 1}
\end{equation} 
and, hence, homothetic $D_n$-invariant curves appear one by one as $\alpha$ is decreased. Moreover, as $\alpha$ approaches $0$, each of theses homothetic curves converges to the regular $n$-gon. Andrews's statement, already mentioned above, that the isoperimetric quotient of a typical curve approaches 0 is not contradicting the existence of the $D_n$-symmetric homothetic solutions because these, according to Andrews, are not attractors.  To illustrate these results we now investigate the evolution of initial curves $\Gamma$ with $D_n$-symmetry for different values of $0 < \alpha < \frac{1}{8}$. Figure \ref{fig:polygonal_kappa01} to Figure \ref{fig:polygonal_kappa001} show evolutions of $D_n$-symmetric curves for $n = 3, 4, 5, 6$, namely for those given by 
\begin{equation}
\label{eq:standardsymmcos}
r(\varphi) = 0.5 + \delta \cos(n \varphi)
\end{equation}
with $n = 3$ and $\delta = 0.05$, with $n = 4$ and $\delta = 0.03$, with $n = 5$ and $\delta = 0.02$, and with $n = 6$ and $\delta = 0.01$. The evolutions were computed with an arclength-uniform mesh with $N=60$ points for $n = 3, 5, 6$ and $N=64$ points for $n = 4$; this ensures that the $2n$ extrema of the $D_n$-symmetric curve given by \ref{eq:standardsymmcos} occur at mesh points in all four cases. Here $S=10^{-4}$ and remeshing was performed at every 10th time step. In each plot we show the first, the last and a selection of intermediate curves with their respective meshes; the computed final shape of the respective evolution is shown separately at the bottom right of each panel. Also shown are the trajectories of the stationary points: they all remain on the respective symmetry axes, thus, demonstrating that symmetry is indeed preserved during the computed evolutions. 

Figure \ref{fig:polygonal_kappa01} illustrates the computation for $\alpha = 0.1$ and shows the evolutions of the above curves for $n = 3$ and $n = 4$. Here $\alpha = \frac{1}{10} < \frac{1}{8} = \frac{1}{3^2-1}$, and the curve with $D_3$-symmetry in panel~(a) becomes more triangular in the sense that the curvature at the three points of maximal curvature increases considerably. The $D_4$-symmetric curve in panel~(b), on the other hand, converges to the circle since $\frac{1}{4^2-1} = \frac{1}{15} < \frac{1}{10} = \alpha$; in particular, the curvature at the four points of maximal curvature decreases.

Figure \ref{fig:polygonal_kappa005} illustrates the computations for $\alpha = 0.05$ and shows the evolutions of the above curves for $n = 3$ up to $n = 6$. Now $\alpha = \frac{1}{20} < \frac{1}{15} = \frac{1}{4^2-1} < \frac{1}{3^2-1} = \frac{1}{8} $, and the curves with $D_3$-symmetry in panel~(a) and with $D_4$-symmetry in panel~(b) approach a triangular and a square shape, respectively. In particular, the final curve with $D_3$-symmetry appears to develop corners, that is, large curvature maxima, and seems already quite close to the equilateral triangle (regular 3-gon). The final curve with $D_4$-symmetry is still further away from the square (regular 4-gon), but the curvature at the four points of maximal curvature is clearly increasing. The curves with $D_5$ and $D_6$ symmetry in panels~(c) and~(d), on the other hand, converge to the circle since $\frac{1}{5^2-1} = \frac{1}{24} < \frac{1}{20} = \alpha$. Note that the convergence to the circle of the $D_5$-symmetric curve is quite slow: while the curve in Figure \ref{fig:polygonal_kappa005}(c) has not quite converged to the circle, the curvature at the five points of maximal curvature is certainly decreasing.

Finally, Figure \ref{fig:polygonal_kappa001} is for $\alpha = 0.01$. Now $\alpha = \frac{1}{100} < \frac{1}{35} = \frac{1}{6^2-1}$ and the evolutions of all of the above curves for $n = 3$ up to $n = 6$ can be seen to approach the respective regular $n$-gon. Here $\alpha$ is quite close to $0$ already and the points of maxima of the curvature are well on their way to becoming corners. This means that the radial distance function $r_{\mbox{\tiny $C$}}$ is getting quite close to losing differentiability at the corresponding $\varphi$-values when the computations stop.

\begin{figure}[t!]
  \centering
  \includegraphics[scale=0.95]{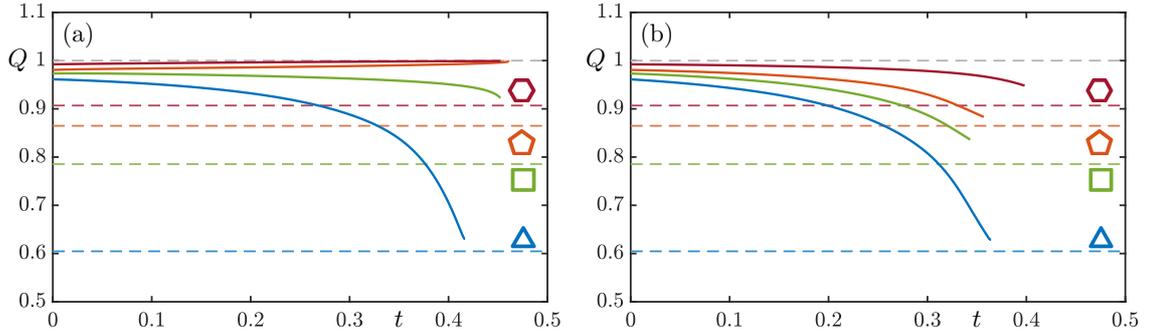}
  \caption{Evolution under \ref{eq:abflow} with $c=0$ and $\alpha = 0.05$ (a) and $\alpha = 0.01$ (b) of the isoperimetric quotient $Q$ for the $D_n$-symmetric curves
from Figure \ref{fig:polygonal_kappa005} and Figure \ref{fig:polygonal_kappa001}, respectively; the horixontal dashed lines are at the values $Q_n$ of the regular $n$-gons.}
  \label{fig:polygonal_isoperimetric}
\end{figure}

Figure \ref{fig:polygonal_kappa01} to Figure \ref{fig:polygonal_kappa001} are numerical evidence that a typical curve $\Gamma \in \mathscr{C}$ with $D_n$-symmetry approaches for $0 < \alpha < \frac{1}{n^2 - 1} \leq \frac{1}{3}$ a curve that is quite close to the regular $n$-gon. As further confirmation of this conclusion, we show in Figure \ref{fig:polygonal_isoperimetric} the time series of the isoperimetric quotient  $Q$ of the $D_n$-symmetric curves for $n = 3$ up to $n = 6$ from Figure \ref{fig:polygonal_kappa005} and Figure \ref{fig:polygonal_kappa001}. In Figure \ref{fig:polygonal_isoperimetric}(a) for $\alpha = 0.05$ the isoperimetric quotient $Q$ converges to 1 when $n = 5$ and $n = 6$, because these curves approach a circle under evolution. The curves for $n = 3$ and $n = 4$ in panel~(a) and all four curves in panel~(b) for $\alpha = 0.01$, on the other hand, approach values $Q^{\star}_n$ in the range $Q_n < Q^{\star}_n <1$ where
\begin{equation}
\label{eq:isongon}
Q_n = \frac{\pi}{n} \cot(\frac{\pi}{n})
\end{equation}
is the isoperimetric quotient of the  regular $n$-gon. The latter values are represented by the horizontal lines in Figure \ref{fig:polygonal_isoperimetric} at
\begin{equation*}
\begin{aligned}
Q_3 & = \frac{\sqrt{3} \pi}{9} \approx 0.6046, \quad 
Q_4 = \frac{\pi}{4} \approx 0.7854, \\
Q_5 & = \frac{\pi}{5} \sqrt{1 + \frac{2}{5}\sqrt{5}} \approx 0.8648 
\quad {\rm and} \quad
Q_6 = \frac{\sqrt{3} \pi}{6} = \approx 0.9069.
\end{aligned}
\end{equation*}
The convergence of $Q$ appears to be faster the smaller $\alpha$; compare panels~(a) and (b).  We observe that, before the computation stops, each computed time series reaches a constant value seemingly above the limit $Q^{\star}_n$, which in turn is above the value $Q_n$ of the $n$-gon shown in Figure \ref{fig:polygonal_isoperimetric}. We found that this plateau-ing results from the finite mesh representation of the curve with piecewise polynomials of a given overall smoothness, in light of large curvature at the $n$ points of maximal curvature.

In accordance with the result of Andrews, we have that $Q^{\star}_n \to Q_n$ as $\alpha \to 0$. In other words, at the $n$ points of its maxima the curvature goes to infinity as $\alpha$ approaches 0. On the other hand, it is also clear that $Q^{\star}_n = Q_n$ cannot be achieved for any $0 < \alpha$ because polygons are not invariant under \ref{eq:abflow}. It remains a considerable challenge to push the calculations further to determine the shape and the associated isoperimetric quotient $Q^{\star}_n$ of the limiting homothetic curve as a function of $\alpha$ satisfying \ref{eq:dncondition}.

\subsubsection{Convergence to $D_n$-symmetric curves}
\label{sec:Cncurves}

\begin{figure}[t!]
  \centering
  \includegraphics[scale=1]{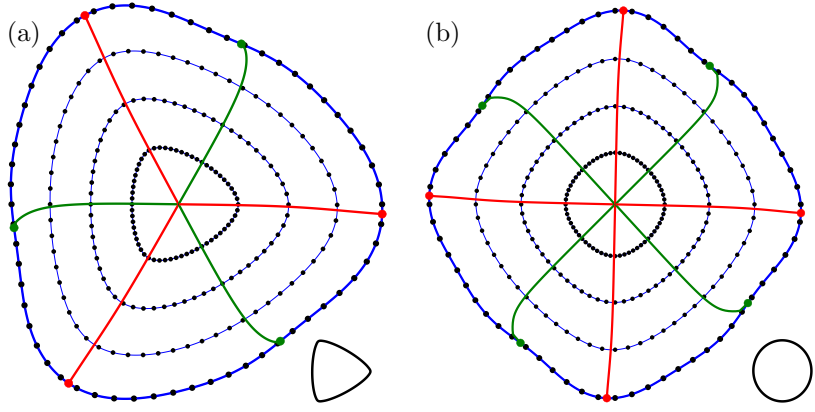}
  \caption{Evolutions under \ref{eq:abflow} with $c=0$ and $\alpha = 0.1$ of the $C_n$-symmetric curves \ref{eq:standardsymmcos_cn} with $n = 3$, $\delta = 0.05$ and $\xi = 0.01$ (a), and with $n = 4$, $\delta = 0.03$ and $\xi = 0.01$ (b), computed with an arclength-uniform mesh with $N=60$ for $n =3$ and $N = 64$ for $n = 4$, remeshing at every 10th time step and $S=10^{-5}$. Shown are a selection of curves with their respective meshes and the trajectories of the stationary points; the computed final shape of the respective evolution is shown at the bottom right of each panel.
}
  \label{fig:cnsymmetry_kappa01}
\end{figure}

Figure \ref{fig:cnsymmetry_kappa01} shows the evolutions under \ref{eq:abflow} with $c=0$ and $\alpha = 0.1$ of two curves that are invariant under  the rotational subgroup $C_n \subset D_n$. More specifically, the curves are defined by 
\begin{equation}
\label{eq:standardsymmcos_cn}
r(\varphi) = 0.5 + \delta \cos(n \varphi) + \xi \cos(2n \varphi + \frac{\pi}{5})
\end{equation}
with $n = 3$ and $\delta = 0.05$ and $\xi = 0.01$ in panel (a), and with $n = 4$, $\delta = 0.03$ and $\xi = 0.01$ in panel (b), respectively. The computations use an arclength-uniform mesh with $N=60$ points for $n = 3$ and $N=64$ points for $n = 4$ with $S=10^{-5}$ and remeshing at every 10th time step. Note that the parameter $\xi$ in \ref{eq:standardsymmcos_cn} breaks the reflectional symmetry of the $D_n$-symmetric curves defined by \ref{eq:standardsymmcos_cn}; hence, the curves shown in Figure \ref{fig:cnsymmetry_kappa01} are symmetry broken perturbations of the curves from Figure \ref{fig:polygonal_kappa01}. During the evolutions shown in Figure \ref{fig:cnsymmetry_kappa01} each curve with $C_n$-symmetry approaches the same limit as the corresponding curve with $D_n$-symmetry. In particular, trajectories of critical points swiftly approach lines through the final point with angles $\frac{\pi}{n}$ between them. This suggests that for $n\geq 3$, if condition \ref{eq:dncondition} on $\alpha$ is satisfied, the subspace of $D_n$-symmetric curves is attracting in the space of $C_n$-symmetric curves. This illustrates Theorem 7.5. from \cite{Andrews2002b} stating that, curves with $C_n$-symmetry converge to the homothetic solutions with $D_n$-symmetry for any $\alpha$ below the bound given by \ref{eq:dncondition}.

Our computations and observations also illustrate a phenomenon related to Bloore's theorem for the Bloore flow \ref{eq:bloore1} on the stability of the sphere in \cite{Bloore1977} where he showed that Fourier terms decay in reverse order. Since each single Fourier term corresponds to a $D_n$-symmetric shape, we expect that $D_n$ symmetry (the symmetry of the last surviving Fourier term) will be attracting in the Bloore flow. What we observed here is in agreement with the results of Andrews \cite{Andrews1996,Andrews2002b,Andrews2002a} and suggests that this property may be more general and apply to a broader class of flows.

\section{Conclusions and Outlook}
\label{sec:conclusions}

We presented an algorithm for computing evolutions of planar star-like curves under a geometric flow, where each curve is represented by a piecewise polynomial periodic radial distance function $r(\varphi)$ with respect to a chosen reference point. Hence, all information of the curve, including the positions of critical points of $r$ and of the curvature, can be computed from this representation explicitly without the need for numerical differentiation. Employing a phase-uniform or an  arclength-uniform mesh, as appropriate, and mesh refinement along increasingly smaller timesteps as determined from the overall curvature, allows us to compute the evolution of the shape of the curve as well as the trajectories of critical points. This was demonstrated with several test-case examples, and two dedicated studies of specific properties of the Andrews-Bloore flow \ref{eq:abflow}. 

More specifically, we illustrated that the number $n_{\mbox{\tiny $C$}}(t)$ of stationary points of a typical curve decreases under the curve-shortening flow \ref{eq:curveshortening} (given by \ref{eq:abflow} with $c=1$ and $\alpha = 1$) to four, where the remaining four stationary points approach the final point tangent to two perpendicular lines. We also confirmed  \ref{conj_U} stating that the ultimate point $U$ produces a monotonic evolution of $n_{\mbox{\tiny $U$}}(t)$ that ends in a manner analogous to the evolution of the stationary points of the curve. Moreover, we presented evidence for \ref{conj_kappa} stating that the number $n_\kappa(t)$ of critical points of the curvature is also monotonically decreasing and, likewise, ends at four as the evolution of the stationary points. Finally, we investigated evolutions of different types of curves under \ref{eq:abflow} for $c=0$ and $0 < \alpha < 1$, which is a case studied extensively by Andrews \cite{Andrews1996,Andrews2002b,Andrews2002a}. Our computations suggest that our conjecture regarding the stationary points, which was formulated for $\alpha=1$, holds more widely over the entire range $0 < \alpha < 1$, including when $0 < \alpha < \frac{1}{3}$ when the limit shape of the curve is not a circle but an interval. This, in turn, can be phrased as the statement that the space of $D_2$-symmetric curves is normally hyperbolic; more generally, we showed that the space of $D_n$-symmetric curves is attracting in the space of $C_n$-symmetric curves. 

There are clearly quite a number of avenues for future research. For example, smooth curves may have limits with developing corners of very large or even infinite curvature, while other parts of the limiting shape may become virtually flat. This poses a problem for the present implementation, as it requires very fine meshes, which are presently chosen to be uniform in phase or arclength. A way of getting closer to such limits is to implement (possibly different types of) curvature-dependent meshes; as was done recently in the context of finite difference methods \cite{Mackenzie2019}. This and other improvements, such as the use of other integrators or interpolating functions, for example, those based on Chebyshev polynomials \cite{chebfunbook,chebfunperiod}, can be formulated and tested within the framework of the general algorithm as presented here. 

A challenging and practically very relevant direction for future research is the development of an algorithm for computing shape and critical point evolution of  compact surfaces under curvature flows. The general framework and the algorithm presented here have been developed with the evolution of surfaces in mind, such as pebbles under abrasion. Indeed, they generalize to this case (and, in fact, to compact surfaces in any dimension) as follows. Any star-like surface can be represented uniquely by a distance function over the sphere $\bbS^2$ with respect to a suitable reference point. The normal and the considered type of curvature at every point can be computed from this spherical representation. The general algorithm then applies in the same way: at each step the surface is represented piecewise by functions of two variables (of some chosen class) over the triangles of a triangulation with respect to a chosen a mesh on $\bbS^2$; the mesh is then evolved by an Euler step at the mesh point; subsequently, the new representation by piecewise polynomials is computed with respect to the reference point (which may be fixed or again the moving centroid), possibly with remeshing. While this general scheme is quite straightforward at a conceptual level, any implementation faces a number of serious challenges. First of all, one needs to choose variables on the sphere, such as Euler angles or stereographic coordinates, as arguments of the distance function and this introduces singularities. A related problem is that generating a uniform mesh on the sphere or more general compact surface without boundary is not as straightforward as for the circle. Finally, the most serious challenge is the choice of the space of functions that are used to represent the surface piecewise over the triangles of a triangulation of the sphere or surface as induced by the chosen mesh. Linear elements give a continuous representation of the surface in the form of a polytope, which may have large numbers of spurious critical points\cite{Domokos2019a}; moreover, the curvature is encoded only by the angle deficit at the mesh point. What is needed is an approximate surface of a prescribed overall smoothness to ensure a more accurate approximation of shape indicators. Requiring differentiability across the edges of the triangulation seems natural, but polynomials of a fixed degree are not able to satisfy this condition (at infinitely many points). One approach may be to require smoothness only at certain points along an edge, or to consider entirely different classes of functions; for example, Chebyshev polynomials over the sphere could be promising in this regard \cite{chebfun3D}. The key challenge is to find a good compromise between choosing a workable, reasonably low-dimensional function space that still generates a faithful representation of shape indicators and, in particular, stationary points and critical points of the curvature.

\section*{Acknowledgements}

EF and GD are grateful to the University of Auckland for their warm reception and support. This research was supported by the Hungarian National Research, Development, and Innovation Office Grant K134199 and the Hungarian National Research, Development, and Innovation Office TKP2020 IE Grant BME Water Sciences \& Disaster Prevention. BK thanks Budapest University of Technology and Economics for their hospitality during several research visits and acknowledges support for his research by Royal Society Te Ap\={a}rangi Marsden Fund grant 16-UOA-286.

\appendix
\section{Arclength parameterization and existence of polar parameterization}
\label{sec:polarmap}

Consider a planar curve $\Gamma$ that is closed and simple (has no points of self-intersection). By general theory, $\Gamma$ has an arclength parameterization 
\begin{equation}
\label{eq:gammaarc}
\begin{array}{rcl}
\tilde{\gamma}: [0,L]  & \to & \bbR^2 \\
                         s \ & \mapsto & \tilde{\gamma}(s) 
= \begin{pmatrix} \tilde{\gamma}_x(s) \\\tilde{\gamma}_y(s) \end{pmatrix} 
\end{array}
\end{equation}
with unit tangent vectors, that is,
\begin{equation}
\label{eq:tangarc}
\|\frac{d\tilde{\gamma}}{ds}(s) \| = 1.
\end{equation}
Here $s$ is the arclength parameter, $L$ is its total arclength and, as before, we assume that $\Gamma$ is sufficiently smooth so that all required derivatives exist. Because $\Gamma$ is closed, $\tilde{\gamma}(0)=\tilde{\gamma}(L)$ and the interval $[0,L]$ is a fundamental interval of the covering space $\bbR$ of the circle $L\bbS^1$. In other words, the parmeterization $\tilde{\gamma}$ is an embedding (because $\Gamma$ is simple) of the circle into the plane, that is, $\Gamma = \tilde{\gamma}(L\bbS^1)$. For the arclength parametrization the normal is
\begin{equation}
\label{eq:narc}
\tilde{\n}(s) = i \frac{d\tilde{\gamma}}{ds}(s) 
\end{equation}
and the (signed) curvature is
\begin{equation}
\label{eq:kappasignarc}
\tilde{\kappa}(s) = \tilde{\n}(s) \cdot \frac{d^2\tilde{\gamma}}{ds^2}
= {\rm det}[ \frac{d\tilde{\gamma}}{ds} \,\, \frac{d^2\tilde{\gamma}}{ds^2}] = \frac{d\tilde{\gamma}_x}{ds} \frac{d^2\tilde{\gamma}_y}{ds^2} - \frac{d\tilde{\gamma}_y}{ds} \frac{d^2\tilde{\gamma}_x}{ds^2}.
\end{equation}
Throughout, the tilde indicates that the respective function is a function of the arclength parameter $s$. Note that the normal is the counter-clockwise rotation of the tangent, and we assume the usual orientation convention for a simple closed curve that the normal $\tilde{\n}(s)$ as defined by \ref{eq:narc} points inwards. Hence, the vector $\frac{d^2\tilde{\gamma}}{ds^2}$ is directed towards the center of curvature, that is, the center of the osculating circle which lies in ${\rm span}(\tilde{\n}(s))$. This means that $\tilde{\kappa}(s)$ is positive when the curve $\Gamma$ is curved inwards and negative when $\Gamma$ is curved outwards. The curvature changes generically from inwards to outwards at (regular) zeros of $\tilde{\kappa}(s)$, which are inflection points of the curve $\Gamma$. 

While the arclength parameterization of a curve $\Gamma$ is ideal from a theory perspective, using it in an algorithm requires one to find it in explicit form. For a simple closed curve this means in practice that one needs to find a parameterization of $\Gamma$ over the circle, which can then be turned into an arclength parameterization by normalization of the periodic variable $\varphi$ with the arclength integral from \ref{eq:polartoarc}. This is achieved by the polar parameterization with respect to a chosen reference point as introduced in \ref{sec:propGamma}. 

We now further clarify when a polar parameterization exists and what points can be chosen as reference points. For any given simple closed curve $\Gamma$ with arclength parameterization \ref{eq:gammaarc} the \textit{polar map} with pole $P \in \bbR^2$, given by
\begin{equation}
\label{eq:lift}
\begin{array}{rcl}
L\, \bbS^1 & \to & (\bbR^+, 2\pi \bbS^1) \\
     s & \mapsto & (\tilde{r}_{\mbox{\tiny $P$}}(s), \tilde{\varphi}_{\mbox{\tiny $P$}}s)) \\[2mm]
&& {\rm with} \ \tilde{r}(s) = \tilde{r}_{\mbox{\tiny $P$}}(s) := \| \tilde{\gamma}(s) - P \| \\[1mm]
&& {\rm and}\ 
\tilde{\varphi}(s) = \tilde{\varphi}_{\mbox{\tiny $P$}}(s) \  {\rm defined \ by} \ e^{i \tilde{\varphi}_{\mbox{\tiny $P$}}(s)} =  \frac{\tilde{\gamma}(s) - P}{\| \tilde{\gamma}(s) - P \|},
\end{array}
\end{equation}
assigns to any point $\tilde{\gamma}(s) \in \Gamma$ the radial radial distance $\tilde{r}(s) = \tilde{r}_{\mbox{\tiny $P$}}(s)$ from and angle $\tilde{\varphi}(s) = \tilde{\varphi}_{\mbox{\tiny $P$}}(s)$ with respect to the pole $P$. Importantly, the polar map \ref{eq:lift} exists and is well defined for any point $P \in \bbR^2$; note that the radial distance function $\tilde{r}$ is strictly positive for $P \not\in \Gamma$.

When $P$ lies outside $\Gamma$ the winding number of $\Gamma$ with respect to $P$ is zero and the graph of the component map $\tilde{\varphi}$ is a $1\!:\!0$ torus knot on the torus $L \bbS^1 \times 2\pi \bbS^1$; in particular, $\tilde{\varphi}$ is not a circle map. For $P$ in its interior, on the other hand, the curve $\Gamma$ winds around $P$ exactly once and $\tilde{\varphi}$ is a circle map (and orientation preserving because of the convention that $\n(s)$ points inward). In particular, its graph is a $1\!:\!1$ torus knot on $L \bbS^1 \times 2\pi \bbS^1$ and  $\tilde{\varphi}$ is surjective. However, the periodic map $\tilde{\varphi}$ is generally not injective because it may have critical points where $\frac{d\tilde{\varphi}}{ds} = 0$, which, generically, come in pairs of a minimum and a maximum; geometrically, at an extremum $s^*$ of $\tilde{\varphi}$ there is a tangency of the ray though $P$ with angle $\tilde{\varphi}(s^*)$ with the curve $\Gamma$. A pair of extrema (dis)appears generically (when the pole $P$ or the curve $\Gamma$ are varied smoothly) at a cubic critical point where also $\frac{d^2\tilde{\varphi}}{ds^2}(s^*) = 0$; compare with \ref{sec:evolution}. 

This discussion shows that the map $\tilde{\varphi}$ is injective, that is, does not have critical points, if and only if all rays from the interior point $P$ intersect the curve $\Gamma$ transversely. This exactly means that $P$ can be chosen as a reference point and, hence, $\Gamma$ is star-like. Since $\tilde{\varphi}$ is surjective, its injectivity means that the inverse $\tilde{\varphi}^{-1}$ exists, so that the arclength $s$ can be expressed as a function of $\varphi \in [-\pi, \pi]$. The polar parameterization with reference point $O=P$ is then 
\begin{equation}
\label{eq:arctopolar}
r_{\mbox{\tiny $O$}}(\varphi) = \tilde{r}_o(\tilde{\varphi}^{-1}(\varphi)) = \tilde{r}_o(s(\varphi))
\end{equation}
with $\gamma_o(\varphi) = r_{\mbox{\tiny $O$}}(\varphi) e^{i\varphi}$ as given by \ref{eq:gammapol}; here $r_{\mbox{\tiny $O$}}$ and $\gamma_o$ appear without a tilde to indicate that they are functions of the phase angle $\varphi$.

\begin{figure}[t!]
  \centering
  \includegraphics[scale=1]{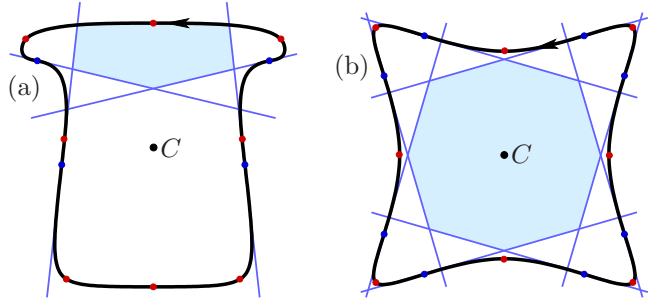}
  \caption{A star-like curve whose centroid $C$ is not a reference point (a), and a nonconvex star-like curve whose centroid $C$ can be chosen as a reference point (b). Also shown are the stationary points (red dots), the inflection points (blue dots) and their tangent lines (blue lines), which bound the set of reference points (shaded).}
  \label{fig:curveexamples}
\end{figure}

The condition that all rays from $O$ intersect the curve transversely is an open property and, hence, all reference points form an open subset of the interior of $\Gamma$. Geometrically, if a point light source is placed at a point reference point then the curve $\Gamma$ is illuminated in its entirety (here we use straight-ray optics). In contrast, a point light source placed at a point from the open complement of the set of reference points will leave parts of $\Gamma$ unlit because there will be at least one ray that has multiple intersections with $\Gamma$ (which hence overshadows itself). This implies that the set of reference points is bounded by the tangent lines through the inflection points of $\Gamma$, where the curvature $\tilde{\kappa}$ is zero and the curve changes from being curved inwards to outwards. In particular, the boundary of the set of reference points is a union of straight line segments and arcs of $\Gamma$ that are positively curved. It follows that the set of reference points is convex; it may or may not contain the centroid $C$ of the curve $\Gamma$, as is shown with two examples in Figure \ref{fig:curveexamples}. A convex curve does not have inflection points and any point in its interior can serve as a reference point for polar coordinates. This shows that star-like curves, which are exactly the curves that admit a polar parameterization, are indeed a mild generalization of convex curves. 

\medskip
\bibliography{references_fdk.bib}

\begin{thebibliography}{10}

\bibitem{ANW67}
J.~Ahlberg, E.~Nilson, and J.~Walsh.
\newblock {\em The Theory of Splines and Their Applications}.
\newblock Elsevier {BV}, 1967.

\bibitem{Alvarez}
L.~Alvarez, F.~Guichard, P.-L. Lions, and J.-M. Morel.
\newblock Axioms and fundamental equations of image processing.
\newblock {\em Archive for Rational Mechanics and Analysis}, 123:199–257,
  1993.

\bibitem{Andrews1994}
B.~Andrews.
\newblock {Contraction of convex hypersurfaces in {E}uclidean space}.
\newblock {\em Calculus of Variations and Partial Differential Equations},
  2(2):151--171, 1994.

\bibitem{Andrews1996}
B.~Andrews.
\newblock {Contraction of convex hypersurfaces by their affine normal}.
\newblock {\em Journal of Differential Geometry}, 43(2):207--230, 1996.

\bibitem{Andrews1998}
B.~Andrews.
\newblock {Evolving convex curves}.
\newblock {\em Calculus of Variations and Partial Differential Equations},
  7(4):315--371, nov 1998.

\bibitem{Andrews2002b}
B.~Andrews.
\newblock {Classification of limiting shapes for isotropic curve flows}.
\newblock {\em Journal of the American Mathematical Society}, 16(02):443--460,
  dec 2002.

\bibitem{Andrews2002a}
B.~Andrews.
\newblock {Non-convergence and instability in the asymptotic behaviour of
  curves evolving by curvature}.
\newblock {\em Communications in Analysis and Geometry}, 10(2):409--449, 2002.

\bibitem{Andrews2002}
B.~Andrews and M.~Feldman.
\newblock {Nonlocal geometric expansion of convex planar curves}.
\newblock {\em Journal of Differential Equations}, 182(2):298--343, 2002.

\bibitem{ArnoldC}
V.~I. Arnold.
\newblock {\em Singularity Theory}.
\newblock Cambridge University Press, 2013.

\bibitem{Barrett2017}
J.~Barrett, K.~Deckelnick, and V.~Styles.
\newblock Numerical analysis for a system coupling curve evolution to reaction
  diffusion on the curve.
\newblock {\em SIAM Journal on Numerical Analysis}, 55(2):1080--1100, 2017.

\bibitem{Barrett2011}
J.~W. Barrett, H.~Garcke, and R.~N{\"u}rnberg.
\newblock The approximation of planar evolution by stable fully implicit finite
  element schemes that equidistribute.
\newblock {\em Numerical Methods for Partial Differential Equations}, 27:1--30,
  2011.

\bibitem{Bloore1977}
F.~J. Bloore.
\newblock The shape of pebbles.
\newblock {\em Journal of the International Association for Mathematical
  Geology}, 9(2):113--122, Apr 1977.

\bibitem{Giblin1}
J.~W. Bruce, P.~J. Giblin, and C.~G. Gibson.
\newblock On caustics of plane curves.
\newblock {\em Amercan Mathematical Monthly}, 88:651--667, 1981.

\bibitem{Bryant1995}
R.~L. Bryant and P.~A. Griffiths.
\newblock Characteristic cohomology of differential systems {II}: Conservation
  laws for a class of parabolic equations.
\newblock {\em Duke Mathematical Journal}, 78(3):531--676, June 1995.

\bibitem{Chow1}
B.~Chow.
\newblock On {H}arnack's inequailty and entropy for the {G}aussian curvature
  flow.
\newblock {\em Communications in Pure and Applied Mathematics}, XLIV:469--483,
  1991.

\bibitem{Daskalopoulos2020}
P.~Daskalopoulos and N.~Sesum.
\newblock Ancient solutions to geometric flows.
\newblock {\em Notices of the American Mathematical Society}, 67(04):1, Apr.
  2020.

\bibitem{Demazure}
M.~Demazure.
\newblock {\em Bifurcations and Catastrophes: Geometry of Solutions to
  Nonlinear Problems}.
\newblock Springer Berlin Heidelberg, 2000.

\bibitem{doedel2007}
E.~J. Doedel.
\newblock Lecture notes on numerical analysis of nonlinear equations.
\newblock In B.~Krauskopf, H.~M. Osinga, and J.~Gal\'an-Vioque, editors, {\em
  Numerical Continuation Methods for Dynamical Systems: Path Following and
  Boundary Value Problems}, pages 1--49. Springer Berlin Heidelberg, 2007.

\bibitem{Domokos2015}
G.~Domokos.
\newblock Monotonicity of spatial critical points evolving under
  curvature-driven flows.
\newblock {\em Journal of Nonlinear Science}, 25(2):247--275, apr 2015.

\bibitem{Domokos2019}
G.~Domokos.
\newblock Natural numbers, natural shapes.
\newblock {\em Axiomathes}, 2019.

\bibitem{DomokosG2012}
G.~Domokos and G.~W. Gibbons.
\newblock The evolution of pebble size and shape in space and time.
\newblock {\em Proceedings of the Royal Society of London A}, 2012.

\bibitem{DomokosLangi2019}
G.~Domokos and Z.~L{\'a}ngi.
\newblock The isoperimetric quotient of a convex body decreases monotonically
  under the {E}ikonal abrasion model.
\newblock {\em Mathematika}, 65(1):119 -- 129, 2019.

\bibitem{Domokos2019a}
G.~Domokos, Z.~L{\'{a}}ngi, and A.~A. Sipos.
\newblock Tracking critical points on evolving curves and surfaces.
\newblock {\em Experimental Mathematics}, pages 1--20, Mar. 2019.

\bibitem{Ruina}
G.~Domokos, J.~Papadopulos, and A.~Ruina.
\newblock Static equilibria of planar, rigid bodies: is there anything new?
\newblock {\em Journal of Elasticity}, 36(1):59--66, 1994.

\bibitem{pebbles}
G.~Domokos, A.~{\'A}. Sipos, T.~Szab{\'o}, and P.~V{\'a}rkonyi.
\newblock Pebbles, shapes and equilibria.
\newblock {\em Mathematical Geosciences}, 42:9, 2010.

\bibitem{chebfunbook}
T.~A. Driscoll, N.~Hale, and L.~N. Trefethen.
\newblock {\em Chebfun Guide}.
\newblock Pafnuty Publications, Oxford, 2014.

\bibitem{Dziuk1999}
G.~Dziuk.
\newblock Discrete anisotropic curve shortening flow.
\newblock {\em SIAM Journal on Numerical Analysis}, 36(6):1808--1830, 1999.

\bibitem{Giblin2}
D.~L. Fidal and P.~J. Giblin.
\newblock Generic one-parameter families of caustics in the plane.
\newblock {\em Mathematical Procedings of the Cambridge Philosophical Society},
  96:425--432, 1984.

\bibitem{Firey1974}
W.~J. Firey.
\newblock Shapes of worn stones.
\newblock {\em Mathematika}, 21(1):1–11, 1974.

\bibitem{Gage1986}
M.~Gage and R.~S. Hamilton.
\newblock The heat equation shrinking convex plane curves.
\newblock {\em Journal of Differential Geometry}, 23(1):69--96, 1986.

\bibitem{Gage1984}
M.~E. Gage.
\newblock Curve shortening makes convex curves circular.
\newblock {\em Inventiones Mathematicae}, 76(2):357--364, 1984.

\bibitem{GolSchaefer}
M.~Golubitsky and D.~Schaeffer.
\newblock {\em Singularities and Groups in Bifurcation Theory, Volume I}.
\newblock Springer Berlin Heidelberg, 1985.

\bibitem{Grayson1987}
M.~A. Grayson.
\newblock The heat equation shrinks embedded plane curves to round points.
\newblock {\em Journal of Differential Geometry}, 26(2):285--314, 1987.

\bibitem{Hamilton}
R.~Hamilton.
\newblock Three-manifolds with positive {R}icci curvature.
\newblock {\em Journal of Differential Geometry}, 17:255--306, 1982.

\bibitem{Huisken1}
G.~Huisken.
\newblock Flow by mean curvature of convex surfaces into spheres.
\newblock {\em Journal of Differential Geometry}, 20:237--266, 1984.

\bibitem{iserles1996a}
A.~Iserles.
\newblock {\em A First Course in the Numerical Analysis of Differential
  Equations}.
\newblock Cambridge University Press, Cambridge New York, 1996.

\bibitem{Ishiwata2017}
T.~Ishiwata and T.~Ohtsuka.
\newblock Evolution of a spiral-shaped polygonal curve by the crystalline
  curvature flow with a pinned tip.
\newblock {\em Discrete {\&} Continuous Dynamical Systems - B}, 22(11):1--35,
  2017.

\bibitem{KPZ}
M.~Kardar, G.~Parisi, and Y.-C. Zhang.
\newblock Dynamic scaling of growing interfaces.
\newblock {\em Physical Review Letters}, 56:889--892, 1986.

\bibitem{Kneser}
A.~Kneser.
\newblock Bemerkungen {\"u}ber die {A}nzahl der {E}xtrema der {K}r{\"u}mmung
  auf geschlossenen {K}urven und {\"u}ber verwandte {F}ragen in einer nicht
  euklidischen {G}eometrie.
\newblock In {\em Festschrift Heinrich Weber}, pages 170--180. Teubner, 1912.

\bibitem{Koenderink}
J.~Koenderink.
\newblock The structure of images.
\newblock {\em Biological Cybernetics}, 50:363--370, 1984.

\bibitem{Mumford}
C.~Lu, Y.~Cao, and D.~Mumford.
\newblock Surface evolution under curvature flows.
\newblock {\em Journal Visual Communication and Imgage Representation},
  13:65--81, 2002.

\bibitem{MacDonald2016}
G.~MacDonald, J.~A. Mackenzie, M.~Nolan, and R.~Insall.
\newblock A computational method for the coupled solution of reaction-diffusion
  equations on evolving domains and manifolds: Application to a model of cell
  migration and chemotaxis.
\newblock {\em Journal of Computational Physics}, 309:207--226, 2016.

\bibitem{Malladi1997}
R.~Malladi and J.~A. Sethian.
\newblock Level set methods for curvature flow, image enchancement, and shape
  recovery in medical images.
\newblock In {\em Visualization and Mathematics}, pages 329--345. Springer
  Berlin Heidelberg, 1997.

\bibitem{Mackenzie2019}
J.~A. Mckenzie, M.~Nolan, C.~F. Rowlatt, and R.~H. Insall.
\newblock An adaptive moving mesh method for forced curve shortening flow.
\newblock {\em SIAM Journal on Scientific Computing}, 41(2):A1170--A1200, 2019.

\bibitem{Mikula1999}
K.~Mikula and D.~\v{S}ev\v{c}ovi\v{c}.
\newblock Solution of nonlinearly curvature driven evolution of plane curves.
\newblock {\em Applied Numerical Mathematics}, 31(2):191 -- 207, 1999.

\bibitem{Mokhtarian}
F.~Mokhtarian, S.~Abbasi, and J.~Kittler.
\newblock Efficient and robust retrieval by shape content through curvature
  scale space.
\newblock In {\em International Workshop on Image Databases and Multimedia
  Search}, pages 35--42, 1996.

\bibitem{Mokhtarian1}
F.~Mokhtarian and R.~Suomela.
\newblock Robust image corner detection through curvature scale space.
\newblock {\em IEEE Transactions on Pattern Analysis and Machine Intelligence},
  20(12):1376--1381, Dec 1998.

\bibitem{Osher2002}
S.~J. Osher and R.~P. Fedkiw.
\newblock {\em Level Set Methods and Dynamic Implicit Surfaces}.
\newblock Springer Berlin Heidelberg, 2002.

\bibitem{Osher1988}
S.~J. Osher and J.~A. Sethian.
\newblock Fronts propagating with curvature-dependent speed: Algorithms based
  on {H}amilton-{J}acobi formulations.
\newblock {\em Journal of Computational Physics}, 79(1):12--49, Nov. 1988.

\bibitem{Perelman}
G.~Perelman.
\newblock {R}icci flow with surgery on three-manifolds.
\newblock {\em http://arXiv.org/math.DG/0303109v1}, 2003.

\bibitem{Popinet1999}
S.~Popinet and S.~Zaleski.
\newblock A front-tracking algorithm for accurate representation of surface
  tension.
\newblock {\em International Journal for Numerical Methods in Fluids},
  30(6):775--793, 1999.

\bibitem{PostonStewart}
T.~Poston and I.~Stewart.
\newblock {\em Catastrophe Theory and its Applications}.
\newblock Dover Publications, Mineola, N.Y, 1996.

\bibitem{rivlin}
T.~Rivlin.
\newblock {\em An Introduction to the Approximation of Functions}.
\newblock Dover, New York, 1981.

\bibitem{Sethian1999}
J.~A. Sethian.
\newblock {\em Level Set Methods and Fast Marching Methods: Evolving Interfaces
  in Computational Geometry, Fluid Mechnics, Computer Vision and Materials
  Science}.
\newblock Cambridge University Press, 1999.

\bibitem{Szaboetal}
T.~Szab{\'o}, G.~Domokos, J.~P. Grotzinger, and D.~J. Jerolmack.
\newblock Reconstructing the transport history of pebbles on {Mars}.
\newblock {\em Nature Communications}, 6:8366, 2015.

\bibitem{chebfun3D}
A.~Townsend, H.~Wilber, and G.~B. Wright.
\newblock Computing with functions in spherical and polar geometries i. the
  sphere.
\newblock {\em SIAM Journal on Scientific Computing}, 2016.

\bibitem{trefethenbook}
L.~N. Trefethen.
\newblock {\em Approximation Theory and Approximation Practice}.
\newblock SIAM, Philadelphia, 2013.

\bibitem{VarkonyiDomokos2006}
P.~V{\'a}rkonyi and G.~Domokos.
\newblock Static equilibria of rigid bodies: Dice, pebbles, and the
  {P}oincar{\'e}-{H}opf theorem.
\newblock {\em Journal of Nonlinear Science}, 16(3):255--281, May 2006.

\bibitem{chebfunperiod}
G.~B. Wright, M.~Javed, H.~Montanelli, and L.~N. Trefethen.
\newblock Extension of {C}hebfun to periodic functions.
\newblock {\em SIAM Journal on Scientific Computing}, 2015.

\end{thebibliography}

\end{document}